\documentclass[11pt,letterpaper,reqno]{amsart}
\usepackage{amsmath,amssymb,amsthm,amsfonts,mathtools}
\usepackage{comment}
\usepackage[margin=1.1in]{geometry}
\usepackage{microtype}
\usepackage{esint}
\makeatletter
\renewcommand{\fint}{%
  \mathchoice
    {\fint@horizontal\displaystyle\textstyle}
    {\fint@horizontal\textstyle\scriptstyle}
    {\fint@horizontal\scriptstyle\scriptscriptstyle}
    {\fint@horizontal\scriptscriptstyle\scriptscriptstyle}%
  \!\int}
\newcommand{\fint@horizontal}[2]{%
  {\setbox\z@=\hbox{$\m@th#1{#2-}{\int}$}%
   \vcenter{\hbox{$\m@th#2-$}}%
   \kern-.5\wd\z@}}
\makeatother
\usepackage[
colorlinks=true,urlcolor=blue,linkcolor=blue,citecolor=magenta]{hyperref}
\usepackage[initials]{amsrefs}

\numberwithin{equation}{section}
\newtheorem{theorem}{Theorem}[section]
\newtheorem{proposition}[theorem]{Proposition}
\newtheorem{lemma}[theorem]{Lemma}
\newtheorem{corollary}[theorem]{Corollary}
\theoremstyle{remark}
\newtheorem{remark}{Remark}[section]
\allowdisplaybreaks[2]
\hypersetup{
 pdftitle={Round-point convergence for the sigma-k flow},
 pdfsubject={Contracting curvature flows and a monotone integral}
}
\title[Entropy and convergence for the
$\sigma_k$-curvature flow]
{Entropy and round-point convergence for the
$\sigma_k$-curvature flow}

\author[T. Luo]{Tianci Luo}
\author[Y. Wei]{Yong Wei}
\address{School of Mathematical Sciences, University of Science and Technology of China, Hefei 230026, P.R. China}
\email{Luo\_tianci@mail.ustc.edu.cn}
\email{yongwei@ustc.edu.cn}

\author[R. Zhang]{Ruijia Zhang}
\address{Department of Mathematics, Sun Yat-sen University, Guangzhou, 510275, China}
\email{zhangrj76@mail.sysu.edu.cn}

\author[R. Zhou]{Rong Zhou}
\address{School of Mathematical Sciences, University of Science and Technology of China, Hefei 230026, P.R. China}
\email{zhourong@mail.ustc.edu.cn}


\begin{document}

\begin{abstract}
We prove that every smooth closed strictly convex hypersurface in
$\mathbb R^{n+1}$ contracts to a round point under the
$\sigma_k$-curvature flow for $2\leq k\leq n-1$, without additional
pinching or symmetry assumptions. After rescaling about the extinction
point to mean width two, the hypersurfaces converge smoothly and
exponentially in normalized time to the unit sphere. The proof is based
on an entropy which is invariant under translations and dilations and
admits a representation by weighted boundary areas of orthogonal
projections. A strengthened Codazzi inequality yields its monotonicity.
An upper bound for the entropy controls the inradius relative to the
mean width, while its dissipation controls the spatial variation of
the normalized speed. We also obtain two-sided curvature bounds at
the natural extinction scale. 
\end{abstract}

\maketitle
\tableofcontents

\section{Introduction}\label{sec:introduction}

A classical problem for contracting curvature flows is to determine
whether closed convex hypersurfaces shrink to round points.
Huisken~\cite{Huisken1984} proved that every smooth closed
strictly convex hypersurface in Euclidean space contracts to a
round point under mean curvature flow. Alternative proofs were
given by Hamilton~\cite{Hamilton1994} and
Andrews~\cite{Andrews1994}.
Chow~\cite{Chow1985} established the corresponding result for
the $n$th root of the Gauss curvature, and later proved
round-point convergence for the square root of the scalar
curvature under an initial curvature pinching assumption
\cite{Chow1987}.
Andrews~\cite{Andrews1994} obtained convergence results for a
large class of homogeneous degree-one curvature flows.
In particular, his later pinching estimates~\cite{Andrews2007}
established round-point convergence for the speeds
$\sigma_k^{1/k}$, $1\leq k\leq n$, without an initial
curvature pinching assumption.

The Gauss curvature flow has been studied extensively,
beginning with Firey's model for the abrasion of convex
bodies~\cite{Firey1974}.
Tso~\cite{Tso1985} proved finite-time contraction to a point,
and Andrews~\cite{Andrews1999} established round-point
convergence for convex surfaces in $\mathbb R^3$.
Entropy estimates for the Gauss curvature flow were developed
by Andrews~\cite{Andrews2000}.
In higher dimensions, Guan--Ni~\cite{GuanNi2017} proved smooth
convergence of the volume-normalized Gauss curvature flow to
a strictly convex self-similar solution.
Andrews--Guan--Ni~\cite{AGN2016} obtained corresponding
convergence results for $K^\alpha$ when
$\alpha>1/(n+2)$.
Together with the classification of self-similar solutions by
Brendle--Choi--Daskalopoulos~\cite{BCD2017}, these results
give round-point convergence for $\alpha>1/(n+2)$.
At the affine exponent $\alpha=1/(n+2)$, the normalized flow
converges to an ellipsoid, as proved by Andrews~\cite{Andrews1996}.

For other curvature speeds of homogeneity greater than one,
convergence results have often required additional assumptions
on the initial hypersurface.
Schulze~\cite{Schulze2006} proved round-point convergence for
powers of the mean curvature under suitable initial curvature
pinching assumptions.
In the appendix, jointly with Schn\"urer, he also proved
convergence without pinching for convex surfaces in
$\mathbb R^3$ moving by $H^\alpha$, $1\leq\alpha\leq5$.
Andrews--McCoy~\cite{AndrewsMcCoy2012} obtained convergence
results for more general curvature speeds under sufficiently
strong initial pinching assumptions.
For powers of the scalar curvature, Alessandroni--Sinestrari
\cite{AS2010} proved round-point convergence for
$\sigma_2^\alpha$, $\alpha>1/2$, under an initial curvature
pinching assumption. They conjectured that this assumption
could be removed.
Li--Wang--Wu~\cite{LiWangWu2021} proved spherical convergence
for the speeds $\sigma_k^\alpha$ on strictly convex axially
symmetric hypersurfaces when
$1/k\leq\alpha\leq c(n,k)$, for some $c(n,k)>1/k$.
For the intermediate speeds $\sigma_k$, $2\leq k\leq n-1$,
round-point convergence from arbitrary smooth closed strictly
convex initial hypersurfaces was not previously known without
additional pinching or symmetry assumptions.

Some related results are available without these assumptions.
Guan--Huang--Liu
\cite{GHL2024}*{Theorem~1.1 and Remark~5.1(2)}
proved finite-time contraction to a point for a class of
curvature flows which includes $\sigma_k^\alpha$ for every
$1\leq k\leq n$ and $\alpha>0$.
Gao--Li--Ma~\cite{GaoLiMa2018}*{Theorem~1.4}
proved that smooth closed strictly convex self-similar
shrinkers for the $\sigma_k^\alpha$-flow are spheres when
$1\leq k<n$ and $\alpha\geq1/k$.
Related rigidity results for high powers of more general
curvature functions were obtained by
Gao--Li--Wang~\cite{GaoLiWang2022}.
However, point extinction and the classification of
self-similar shrinkers do not by themselves imply convergence
of the rescaled flow. For this, one needs uniform estimates
and a way to show that the normalized hypersurfaces approach
a sphere.

The curvature pinching arguments used for degree-one flows
do not extend directly to the speeds $\sigma_k$,
$2\leq k\leq n-1$.
As observed in~\cite{AS2010}, the higher homogeneity
introduces gradient terms which are difficult to control
without an initial pinching assumption.
The entropy method for the Gauss curvature flow also uses
features particular to Gauss curvature.
In Gauss coordinates, $K\,d\mu=d\sigma$, whereas
$\sigma_k\,d\mu$ depends on the principal radii of curvature
when $k<n$. A different monotone quantity is therefore
needed for the intermediate flows.

Related results have been obtained for constrained curvature
flows. Andrews--Wei~\cite{AndrewsWei2021} proved smooth
convergence to spheres for flows by positive powers of
$\sigma_k$, preserving the enclosed volume or a mixed volume,
without initial curvature pinching assumptions.
Their proof uses the monotonicity of an isoperimetric ratio
to control the inradius and diameter, followed by a lower
speed estimate.
For further results on contracting convex hypersurfaces,
we refer to~\cite{AMZ2013} and the
monograph~\cite{ACGL2020}.

In this paper, we prove round-point convergence for the
$\sigma_k$-curvature flow, $2\leq k\leq n-1$, starting from arbitrary
smooth closed strictly convex hypersurfaces.
The main ingredient is a strengthened Codazzi inequality
which yields a scale-invariant monotone entropy.
The entropy controls the inradius relative to the mean width,
and its dissipation gives convergence of the normalized
support function before uniform curvature bounds are obtained.
These estimates lead to smooth exponential convergence to
a sphere. 

\subsection{Main result}

Let $X_0:M^n\to\mathbb R^{n+1}$ be a smooth embedding of a
closed connected manifold, with strictly convex image $M_0$.
For the principal curvatures
$\kappa=(\kappa_1,\ldots,\kappa_n)$ of a hypersurface, set
\begin{equation}\label{eq:intro-symmetric}
 \sigma_k(\kappa)
 =\sum_{1\leq i_1<\cdots<i_k\leq n}
       \kappa_{i_1}\cdots\kappa_{i_k},
 \qquad
 E_k(\kappa)=\binom nk^{-1}\sigma_k(\kappa).
\end{equation}
In particular, $\sigma_1=H$ and $\sigma_n=K$, where
$H$ and $K$ are the mean curvature and Gauss curvature,
respectively. Thus $E_1=H/n$ and $E_n=K$.

We consider the contracting $\sigma_k$-curvature flow
of hypersurfaces $M_t=X(M,t)$, written with the
normalized speed $E_k$ as
\begin{equation}\label{eq:flow}
 \begin{cases}
  \partial_tX=-E_k(W)\nu,\\
  X(\cdot,0)=X_0,
 \end{cases}
\end{equation}
where $X(\cdot,t)$ is a family of smooth embeddings,
$\nu$ is the outward unit normal, and $W$ is the shape
operator of $M_t$.
Since $E_k=\binom nk^{-1}\sigma_k$, the two speeds differ
only by a constant rescaling of time.
Let $\Omega_t$ be the convex body enclosed by $M_t$,
and let $h(\cdot,t)$ be its support function.

\begin{theorem}\label{thm:main}
Let $n\geq 3$ and $2\leq k\leq n-1$.
Suppose that $M_0\subset\mathbb R^{n+1}$ is smooth, closed,
connected, embedded, and strictly convex.
Then \eqref{eq:flow} has a unique maximal smooth strictly
convex solution on $[0,T)$, where $T<\infty$, and
$M_t$ contracts to a point $x_*\in\mathbb R^{n+1}$.

Let $b(t)$ be half the mean width of $\Omega_t$, and define
\begin{equation}\label{eq:mean-width-normalization}
\begin{gathered}
 b(t)=\frac1{|\mathbb S^n|}
      \int_{\mathbb S^n}h(z,t)\,d\sigma,
 \qquad
 \tau(t)=\int_0^t b(s)^{-(k+1)}\,ds,\\
 Y(\cdot,\tau(t))=\frac{X(\cdot,t)-x_*}{b(t)}.
\end{gathered}
\end{equation}
Then $\tau(t)\to\infty$ as $t\to T$, and the support
function $u$ of $Y$ has mean one and satisfies
\begin{equation}\label{eq:main-convergence}
 u(\cdot,\tau)\longrightarrow1
 \quad\text{in }C^j(\mathbb S^n)
 \quad\text{for every integer }j\geq0.
\end{equation}
The convergence is exponential in the normalized time.

There are $c,C>0$, depending only on $M_0,n,k$, such that
\begin{equation}\label{eq:main-curvature}
 c(T-t)^{-1/(k+1)}
 \leq\kappa_i(M_t)
 \leq C(T-t)^{-1/(k+1)},
 \qquad 1\leq i\leq n,
\end{equation}
for all $t$ sufficiently close to $T$.
In particular, the flow satisfies the corresponding
Type I upper bounds for $H$ and $E_k$.
\end{theorem}

The existence and point-extinction assertions follow from
Proposition~\ref{prop:point-extinction}. The normalization fixes the mean width at two. The exponential convergence rate is
given in Corollary~\ref{cor:exponential-convergence}.

Theorem~\ref{thm:main} establishes round-point convergence
for the flow~\eqref{eq:flow} in the intermediate range
$2\leq k\leq n-1$, without any initial curvature pinching
or symmetry assumption.
The monotone entropy also provides an alternative proof
of the rigidity of self-similar shrinkers
(Corollary~\ref{cor:shrinker-rigidity})
and yields a quantitative stability estimate
(Proposition~\ref{prop:monotone-stability}).

\subsection{Discussion of the proof}

The main difficulty for the intermediate flows
$2\leq k\leq n-1$ is to obtain convergence without an initial
curvature pinching assumption. Our proof uses a new
monotone entropy to control the geometry of the evolving
hypersurfaces. Its dissipation gives convergence of the
normalized support function before uniform curvature
bounds are established.

\medskip
\noindent\emph{A strengthened Codazzi inequality.}
Let $\sigma$ and $\bar\nabla$ denote the round metric
and connection on $\mathbb S^n$.
The principal-radius tensors of the original and
normalized hypersurfaces are
$B=\bar\nabla^2h+h\sigma$ and
$\widetilde B=\bar\nabla^2u+u\sigma$, respectively.
Set $G(B)=-E_k(B^{-1})=-F$. 
The Codazzi equation implies that $\bar\nabla B$ is
fully symmetric. For $1\leq k<n$,
Proposition~\ref{prop:strict-codazzi} proves
\begin{equation*}
 -\ddot G^{ij,rs}T_{mij}T_{mrs}
 \geq \frac{2}{nE_1(B)}\dot G^{ij}T_{ipp}T_{jqq}
 +\frac1F\sum_m
       \bigl(\dot G^{ij}T_{mij}\bigr)^2
\end{equation*}
for every fully symmetric three-tensor $T$.
If the last term is omitted, equality holds only
for $T=0$, and the coefficient $2$ is optimal. 

Taking $T=\bar\nabla B$, the remainder
$\mathcal C_G$ defined in \eqref{eq:strict-remainder} satisfies
\begin{equation*}
 \mathcal C_G\geq
 \frac{|\bar\nabla F|_\sigma^2}{F},
 \qquad
 \mathcal C_G=0\quad\Longleftrightarrow\quad
 \bar\nabla B=0.
\end{equation*}
Thus the inequality controls the gradient of the speed
and characterizes the equality case.

\medskip
\noindent\emph{A monotone entropy.}
For $2\leq k\leq n-1$, set
\begin{equation}\label{eq:exponents}
 p=\frac{n-1}{k-1},
 \qquad
 q=p-1=\frac{n-k}{k-1},
\end{equation}
and define
\begin{equation}\label{eq:intro-integral}
 \mathcal E(M)
 =\frac1{|\mathbb S^n|}
   \int_M E_k(W)^p E_1(B)^q\,d\mu.
\end{equation}
This entropy is invariant under translations and
dilations. It has a geometric interpretation in
terms of orthogonal projections. Indeed, writing $j=n-k$, we prove
\begin{equation}\label{eq:intro-projection-integral}
 \mathcal E(\partial\Omega)
 =\frac1{|\mathbb S^j|}
  \int_{\operatorname{Gr}(n+1,j+1)}
  \int_{\partial(\Omega|\Pi)}
  \bigl(FE_1(B)\bigr)^q(\nu_\Pi)
  \,d\mu_\Pi\,d\Pi.
\end{equation}
Here $\Omega|\Pi$ is the orthogonal projection of
$\Omega$ onto $\Pi$, and $d\Pi$ is the
rotation-invariant probability measure on
$\operatorname{Gr}(n+1,j+1)$.
The outward unit normal and boundary area element
of the projection are denoted by $\nu_\Pi$ and
$d\mu_\Pi$. The weight is evaluated on the original hypersurface
at the normal $\nu_\Pi$.
Since the projected boundary area has homogeneity
$j$ and $FE_1(B)$ has homogeneity $1-k$,
scale invariance determines
$q=j/(k-1)$ and $p=q+1$.

The Newton--Maclaurin inequalities and the
strengthened Codazzi inequality give
\begin{equation*}
 \mathcal E(M)\geq1,
 \qquad
 \frac{d}{dt}\mathcal E(M_t)\leq0,
\end{equation*}
with equality in either statement only for round
spheres. The entropy also controls the shape of a convex
body. If $b$ is the mean support function of
$\Omega$, Proposition~\ref{prop:noncollapse} gives
\begin{equation*}
 r_-(\Omega)\geq
 c(n,k)b\,\mathcal E(\partial\Omega)^{-(k-1)},
 \qquad
 \operatorname{diam}\Omega\leq C(n)b.
\end{equation*}
Thus an upper bound for the entropy prevents
collapse after fixing the mean width.
The inradius estimate follows from
Cauchy's projection formula and H\"older's inequality. For earlier entropy estimates for curvature flows, see Andrews~\cites{AndrewsEntropy1994,Andrews2000}
and Hamilton~\cite{HamiltonEntropy1994}.

We also obtain a quantitative stability estimate.
If $c$ is the Steiner point of $\Omega$, then
Proposition~\ref{prop:monotone-stability} gives
\begin{equation*}
 \frac{d_H(\Omega,\overline B_b(c))}{b}
 \leq C(n,k)
 \bigl(\mathcal E(\partial\Omega)-1
 \bigr)^{2/(n+2)}.
\end{equation*}
This is an independent geometric consequence of
the entropy and is not needed for the convergence
argument.

\medskip
\noindent\emph{Convergence of the normalized flow.}
We normalize the mean support function to one
as in \eqref{eq:mean-width-normalization}.
The entropy bound gives uniform control of the
inradius and diameter. Together with the Harnack
estimate and a maximum principle, this yields
uniform upper bounds for the normalized speed
and its product with the mean principal radius.

Write $F=E_k(\widetilde B^{-1})$ and
$\overline F=\fint_{\mathbb S^n}F\,d\sigma$.
For the quantity $\mathcal D$ defined in
\eqref{eq:integrated-remainder},
the entropy formula and the strengthened
Codazzi inequality give
\begin{equation*}
 \frac{d}{d\tau}\mathcal E(Y(\tau))
 \leq-\mathcal D(\tau),
 \qquad
 \mathcal D(\tau)\geq
 c\fint_{\mathbb S^n}|\bar\nabla F|^2\,d\sigma
 \geq c\fint_{\mathbb S^n}
              |F-\overline F|^2\,d\sigma.
\end{equation*}
On the other hand, the normalized support function
satisfies
\begin{equation*}
 \partial_\tau(u-1)
 =\overline F(u-1)-(F-\overline F),
 \qquad \overline F\geq1.
\end{equation*}
Integrating this equation towards future times
gives
\begin{equation*}
 \fint_{\mathbb S^n}|u(z,\tau)-1|^2\,d\sigma
 \leq C\int_\tau^\infty
     e^{-(s-\tau)}\mathcal D(s)\,ds
 \longrightarrow0,
\end{equation*}
since
\begin{equation*}
 \int_0^\infty\mathcal D(s)\,ds
 \leq\mathcal E(M_0)-1.
\end{equation*}
Convexity and the diameter bound give
$|\bar\nabla u|\leq C$, so $u$ converges
uniformly to one.

This proves convergence of the support function
without first obtaining two-sided curvature bounds.
In particular, $u\geq1/2$ for all sufficiently
large normalized times. A maximum principle applied
to $ \log(Fu^{k+1})$  
then gives a positive lower bound for $F$.
Combining this with the previous upper bounds yields
\begin{equation*}
 c\sigma\leq\widetilde B\leq C\sigma.
\end{equation*}
Parabolic regularity now gives smooth convergence
to the unit sphere. Linearization at the sphere,
together with control of the translation modes,
gives exponential convergence as stated in
Corollary~\ref{cor:exponential-convergence}.
Rescaling yields the curvature estimates in
Theorem~\ref{thm:main}.

The paper is organized as follows.
Section~\ref{sec:extinction} collects the preliminary
results. Section~\ref{sec:codazzi} proves the
Codazzi inequality.
Section~\ref{sec:integral} establishes the entropy
and its geometric estimates.
Section~\ref{sec:convergence} completes the
convergence proof.

\subsection{AI Declaration}  
During the preparation of this paper, the authors used OpenAI's
ChatGPT (6 Astra) to assist with exploring proof strategies, searching the literature, developing and checking mathematical arguments and calculations, and improving the exposition. The authors critically evaluated the suggestions, checked the mathematical arguments and cited sources, and revised the manuscript accordingly. The authors take full responsibility for all mathematical statements, proofs, references, and conclusions in the final manuscript.

\section{Preliminaries}\label{sec:extinction}

Throughout this section, $n\geq2$ and $1\leq k\leq n$.
We recall the algebraic identities, geometric formulas, and estimates
used below. Repeated tensor indices are summed from $1$ to $n$ unless
otherwise stated.

\subsection{Elementary symmetric functions}\label{subsec:symmetric}

For $\lambda\in\mathbb R^n$ and $0\leq j\leq n$, let
\begin{equation*}
 \sigma_j(\lambda)
 =\sum_{1\leq i_1<\cdots<i_j\leq n}
       \lambda_{i_1}\cdots\lambda_{i_j},
 \qquad E_j=\frac{\sigma_j}{\binom nj}.
\end{equation*}
We set $\sigma_0=E_0=1$ and $\sigma_j=E_j=0$ for $j<0$ or $j>n$.
When fewer variables are specified, $E_j$ is normalized in that number
of variables. We work on the positive cone
\begin{equation*}
 \Gamma_+=\{\lambda\in\mathbb R^n:\lambda_i>0\text{ for all }i\}.
\end{equation*}
Writing $\lambda|i$ and $\lambda|i,l$ for the vectors with the indicated
entries omitted, we have
\begin{equation}\label{eq:symmetric-derivatives}
 \frac{\partial\sigma_j}{\partial\lambda_i}
 =\sigma_{j-1}(\lambda|i),
 \qquad
 \frac{\partial^2\sigma_j}{\partial\lambda_i\partial\lambda_l}
 =\begin{cases}
   \sigma_{j-2}(\lambda|i,l),&i\ne l,\\
   0,&i=l.
  \end{cases}
\end{equation}
In particular,
\begin{equation}\label{eq:symmetric-trace-identities}
 \sum_i\lambda_i\frac{\partial\sigma_j}{\partial\lambda_i}
 =j\sigma_j,
 \qquad
 \sum_i\lambda_i^2\frac{\partial\sigma_j}{\partial\lambda_i}
 =\sigma_1\sigma_j-(j+1)\sigma_{j+1}.
\end{equation}
For a positive definite symmetric matrix $B$, we also have
\begin{equation}\label{eq:inverse-symmetric}
 E_j(B^{-1})=\frac{E_{n-j}(B)}{\det B},
 \qquad 0\leq j\leq n.
\end{equation}

A smooth symmetric function $\Phi$ on $\Gamma_+$ is also regarded as
a function of positive definite symmetric matrices by evaluation at
their eigenvalues. Its matrix derivatives are defined by
\begin{equation*}
 \left.\frac{d}{ds}\Phi(B+sV)\right|_{s=0}
 =\dot\Phi^{ij}V_{ij},
 \qquad
 \left.\frac{d^2}{ds^2}\Phi(B+sV)\right|_{s=0}
 =\ddot\Phi^{ij,pq}V_{ij}V_{pq}.
\end{equation*}
Write $\Phi_i$ and $\Phi_{ij}$ for derivatives in the eigenvalue
variables. At $B=\operatorname{diag}(r_1,\ldots,r_n)$, the spectral
derivative formula \cite{Andrews2007}*{Theorem~5.1} reads
\begin{equation}\label{eq:spectral-hessian}
 \dot\Phi^{ij}=\Phi_i\delta_{ij},
 \qquad
 \ddot\Phi^{ij,pq}V_{ij}V_{pq}
 =\sum_{i,j}\Phi_{ij}V_{ii}V_{jj}
  +2\sum_{i<j}\frac{\Phi_i-\Phi_j}{r_i-r_j}V_{ij}^2.
\end{equation}
The divided differences extend continuously to repeated eigenvalues,
with value $\Phi_{ii}-\Phi_{ij}$ when $r_i=r_j$.

\subsection{Convex hypersurfaces and evolution equations}
\label{subsec:geometry-evolution}

Let $M=\partial\Omega\subset\mathbb R^{n+1}$ be smooth, closed,
and strictly convex, with outward unit normal $\nu$.
Its induced metric, connection, second fundamental form, and area
measure are denoted by $g$, $\nabla$, $h_{ij}$, and $d\mu$.
With $D$ the Euclidean connection, our convention is
\begin{equation*}
 h_{ij}=\langle D_{e_i}\nu,e_j\rangle,
 \qquad W=(h^i{}_j),\quad h^i{}_j=g^{il}h_{lj},
 \qquad H=\operatorname{tr}W.
\end{equation*}
The eigenvalues $\kappa_i$ of $W$ are positive, and
\begin{equation}\label{eq:intrinsic-notation}
 F=E_k(W)=E_k(\kappa).
\end{equation}
Let $\sigma$, $\bar\nabla$, and $d\sigma$ denote the round metric,
connection, and volume element on $\mathbb S^n$, and write
\begin{equation*}
 \fint_{\mathbb S^n}v\,d\sigma
 =\frac1{|\mathbb S^n|}\int_{\mathbb S^n}v\,d\sigma.
\end{equation*}
The support function and inverse Gauss map are
\begin{equation}\label{eq:inverse-gauss}
 h(z)=\sup_{x\in\Omega}\langle x,z\rangle,
 \qquad X(z)=\bar\nabla h(z)+h(z)z.
\end{equation}
Identifying both tangent spaces with $z^\perp$, the principal-radius
tensor is
\begin{equation}\label{eq:radius-tensor}
 B_{ij}=\bar\nabla_i\bar\nabla_jh+h\sigma_{ij},
 \qquad W=B^{-1},
 \qquad E_1(B)=\frac1n\operatorname{tr}_\sigma B.
\end{equation}
Its eigenvalues $r_i=\kappa_i^{-1}$ are the principal radii.
All determinants and contractions involving $B$ use $\sigma$.
In particular,
\begin{equation}\label{eq:radius-identities}
 d\mu=\det B\,d\sigma,
 \qquad E_1(B)=\frac{E_{n-1}(W)}{E_n(W)},
 \qquad F=\frac{E_{n-k}(B)}{\det B}.
\end{equation}
The Codazzi equation
\begin{equation}\label{eq:codazzi-radius}
 \bar\nabla_iB_{jl}=\bar\nabla_jB_{il}
\end{equation}
makes $\bar\nabla B$ fully symmetric.
The associated Newton tensors are divergence free:
\begin{equation}\label{eq:radius-newton-divergence}
 \bar\nabla_i\left(\frac{\partial E_j(B)}{\partial B_{il}}\right)=0
 \qquad(1\leq j\leq n).
\end{equation}
This includes the cofactor tensor $(\det B)B^{-1}$.

In the normal parametrization $\partial_tX=-F\nu$, the metric
variation is
\begin{equation}\label{eq:area-evolution}
 \partial_tg_{ij}=-2Fh_{ij}.
\end{equation}
See \cite{Huisken1984}*{Section~3} and
\cite{AMZ2013}*{Lemma~9}.
In Gauss coordinates, the support equation is
\begin{equation}\label{eq:support-original}
 \partial_th=G(B),
 \qquad G(B)=-E_k(B^{-1})=-F.
\end{equation}
We write
\begin{equation}\label{eq:coordinate-derivatives}
 \dot G^{ij}=\frac{\partial G}{\partial B_{ij}},
 \qquad
 \ddot G^{ij,pq}=\frac{\partial^2G}{\partial B_{ij}\partial B_{pq}}.
\end{equation}
These derivatives are evaluated at the radius tensor under
consideration, with $\sigma$ fixed. 
Since $E_k^{1/k}$ is inverse-concave
(see \cite{AMZ2013}*{Section~2}),
the function $E_k(B^{-1})^{-1/k}$ is concave.
As $s\mapsto-s^{-k}$ is increasing and concave for $s>0$,
it follows that $G(B)=-E_k(B^{-1})$ is concave.
Moreover, $(\dot G^{ij})$ is positive definite, and
homogeneity gives $\dot G^{ij}B_{ij}=kF$. 
At fixed Gauss normal,
\begin{equation}\label{eq:gauss-velocity}
 \partial_tX=-Fz-\bar\nabla F.
\end{equation}
Since $dX=B$, the tangential velocity relative to the normal
parametrization is $-\bar\nabla F=-B\nabla F$.

\begin{lemma}\label{lem:gauss-evolutions}
Along the original flow \eqref{eq:flow}, with time derivatives taken at fixed Gauss normal,
\begin{align}
 (\partial_t-\dot G^{ij}\bar\nabla_i\bar\nabla_j)h
 &=-(k+1)F+(\dot G^{ij}\sigma_{ij})h,
 \label{eq:evolution-h}\\
 (\partial_t-\dot G^{ij}\bar\nabla_i\bar\nabla_j)F
 &=(\dot G^{ij}\sigma_{ij})F,
 \label{eq:evolution-F}\\
 (\partial_t-\dot G^{pq}\bar\nabla_p\bar\nabla_q)B_{ij}
 &=\ddot G^{pq,rs}\bar\nabla_iB_{pq}\bar\nabla_jB_{rs}
   +(k-1)F\sigma_{ij}
   -(\dot G^{pq}\sigma_{pq})B_{ij}.
 \label{eq:evolution-B}
\end{align}
Define
\begin{equation}\label{eq:strict-remainder}
 \mathcal C_G
 =-\ddot G^{ij,pq}\bar\nabla_mB_{ij}\bar\nabla_mB_{pq}
  -\frac{2n}{E_1(B)}\dot G^{ij}
   \bar\nabla_iE_1(B)\bar\nabla_jE_1(B).
\end{equation}
Then
\begin{equation}
    \begin{aligned}\label{eq:mean-radius-evolution}
 (\partial_t-\dot G^{ij}\bar\nabla_i\bar\nabla_j)E_1(B)
 =&-\frac{2}{E_1(B)}\dot G^{ij}\bar\nabla_iE_1(B)\bar\nabla_jE_1(B)\\
       &+(k-1)F-(\dot G^{ij}\sigma_{ij})E_1(B)
       -\frac1n\mathcal C_G,
 \end{aligned}
\end{equation}
and
\begin{equation}\label{eq:product-original}
    (\partial_t-\dot G^{ij}\bar\nabla_i\bar\nabla_j)
       \bigl(FE_1(B)\bigr)
 =-2\dot G^{ij}\bar\nabla_i\log E_1(B)
\bar\nabla_j\bigl(FE_1(B)\bigr)+(k-1)F^2-\frac Fn\mathcal C_G.
\end{equation}
\end{lemma}

\begin{proof}
The support equation and homogeneity give
\eqref{eq:evolution-h} and \eqref{eq:evolution-F}.
Differentiating \eqref{eq:radius-tensor} gives
$\partial_tB_{ij}=-\bar\nabla_i\bar\nabla_jF-F\sigma_{ij}$.
Commuting derivatives on $\mathbb S^n$ and using
\eqref{eq:codazzi-radius}, we obtain
\begin{equation*}
 \dot G^{pq}\bar\nabla_i\bar\nabla_jB_{pq}
 =\dot G^{pq}\bar\nabla_p\bar\nabla_qB_{ij}
  +kF\sigma_{ij}-(\dot G^{pq}\sigma_{pq})B_{ij}.
\end{equation*}
Together with
\begin{equation*}
 \bar\nabla_i\bar\nabla_jF
 =-\dot G^{pq}\bar\nabla_i\bar\nabla_jB_{pq}
  -\ddot G^{pq,rs}\bar\nabla_iB_{pq}\bar\nabla_jB_{rs},
\end{equation*}
this proves \eqref{eq:evolution-B}. Compare
\cite{AMZ2013}*{Lemma~10}.
Taking the trace gives \eqref{eq:mean-radius-evolution}, and
the product rule gives \eqref{eq:product-original}.
\end{proof}

\subsection{Point extinction and the Harnack estimate}
\label{subsec:extinction-estimates}

We recall the point-extinction theorem of Guan--Huang--Liu
\cite{GHL2024}*{Theorem~1.1 and Remark~5.1(2)}.

\begin{proposition}\label{prop:point-extinction}
Let $M_0\subset\mathbb R^{n+1}$ be smooth, closed, connected,
embedded, and strictly convex. The flow \eqref{eq:flow} has a unique
maximal smooth strictly convex solution on $[0,T)$, where $T<\infty$,
and $M_t$ contracts to a single point $x_*\in\mathbb R^{n+1}$.
Moreover,
\begin{equation}\label{eq:lifetime}
 \frac{r_-(\Omega_0)^{k+1}}{k+1}
 \leq T\leq
 \frac{r_+(\Omega_0)^{k+1}}{k+1},
\end{equation}
where $r_-$ and $r_+$ are the inradius and circumradius.
\end{proposition}

We use the differential Harnack inequality of
Andrews~\cite{AndrewsHarnack1994}*{Theorem~5.6},
in the form recorded in
\cite{AMZ2013}*{Theorem~13}.
The estimate is written at fixed Gauss normal and applies directly
to the speed $E_k$.

\begin{lemma}\label{lem:harnack}
Let $M_t$, $0\leq t<T$, be the maximal smooth strictly convex
solution of \eqref{eq:flow}. Then
\begin{equation}\label{eq:harnack}
 \partial_tF(z,t)+\frac{k}{k+1}\frac{F(z,t)}t\geq0
 \qquad(z\in\mathbb S^n,\ 0<t<T).
\end{equation}
Consequently,
\begin{equation}\label{eq:harnack-integrated}
 F(z,t_2)\geq
 \left(\frac{t_1}{t_2}\right)^{k/(k+1)}F(z,t_1)
 \qquad(0<t_1\leq t_2<T).
\end{equation}
\end{lemma}

\subsection{Mean-width normalization}\label{subsec:normalization}

Let $x_*$ and $T$ be given by Proposition~\ref{prop:point-extinction}.
The mean support $b(t)$ in \eqref{eq:mean-width-normalization}
is unchanged by translation.
We use the normalized flow $Y=(X-x_*)/b(t)$ and write
\begin{align}
 u(z,\tau(t))
 &=\frac{h(z,t)-\langle x_*,z\rangle}{b(t)},
 \label{eq:normalized-support-definition}\\
 \widetilde B_{ij}
 &=\bar\nabla_i\bar\nabla_ju+u\sigma_{ij},
 \qquad \widetilde W=\widetilde B^{-1}.
 \label{eq:normalized-notation}
\end{align}
The enclosed body is $\widetilde\Omega_\tau=(\Omega_t-x_*)/b(t)$.
Then
\begin{equation}\label{eq:normalized-mean-support}
 \fint_{\mathbb S^n}u\,d\sigma
 =\fint_{\mathbb S^n}E_1(\widetilde B)\,d\sigma=1.
\end{equation}
Thus the mean width is two. Equivalently, the mean perimeter of
the planar projections is $2\pi$, since
\begin{equation*}
 \frac1{2\pi}\int_{\operatorname{Gr}(n+1,2)}
 \operatorname{Length}\bigl(\partial(\widetilde\Omega_\tau|\Pi)\bigr)
 \,d\Pi=\fint_{\mathbb S^n}u\,d\sigma.
\end{equation*}
Here $d\Pi$ is the rotation-invariant probability measure. The
identity follows by averaging the planar support-function formula
for perimeter.

Set
\begin{equation}\label{eq:average-normalized-speed}
 \overline F(\tau)
 =\fint_{\mathbb S^n}E_k(\widetilde B^{-1})\,d\sigma.
\end{equation}
Writing $\widetilde r_i$ for the eigenvalues of $\widetilde B$,
the scaling relations are
\begin{equation}\label{eq:scaling-identities}
\begin{gathered}
 B=b\widetilde B,\qquad W=b^{-1}\widetilde W,
 \qquad r_i=b\widetilde r_i,\qquad
 F=b^{-k}E_k(\widetilde B^{-1}),\\
 \frac{dt}{d\tau}=b^{k+1},\qquad
 \frac{db}{dt}=-b^{-k}\overline F,
 \qquad \frac{db}{d\tau}=-\overline F b.
\end{gathered}
\end{equation}
In the last identity $b$ is composed with $t(\tau)$.
The Newton--Maclaurin and Jensen inequalities give
\begin{equation}\label{eq:average-speed-lower}
 \overline F
 \geq\fint_{\mathbb S^n}E_1(\widetilde B)^{-k}\,d\sigma
 \geq\left(\fint_{\mathbb S^n}E_1(\widetilde B)\,d\sigma\right)^{-k}
 =1.
\end{equation}
Lemma \ref{lem:normalized-speed} will show that $\tau(t)\to\infty$ as $t\to T$.

We retain the extinction-time scale only for sphere comparison.
A sphere evolving by \eqref{eq:flow} satisfies
\begin{equation}\label{eq:sphere-ode}
 r'(t)=-r(t)^{-k},\qquad
 r(t)^{k+1}=r(0)^{k+1}-(k+1)t.
\end{equation}
Set
\begin{equation}\label{eq:extinction-scale}
 \ell(t)=\bigl((k+1)(T-t)\bigr)^{1/(k+1)}.
\end{equation}
Comparison with spheres gives
\begin{equation}\label{eq:extinction-radius-comparison}
 r_-(\Omega_t)\leq\ell(t)\leq\operatorname{diam}\Omega_t.
\end{equation}

\begin{lemma}\label{lem:normalized-evolutions}
The normalized support function satisfies
\begin{equation}\label{eq:normalized-support}
 \partial_\tau u
 =\overline F u+G(\widetilde B)
 =\overline F u-E_k(\widetilde B^{-1}).
\end{equation}
Let $\widetilde{\mathcal C_G}$ denote
\eqref{eq:strict-remainder} evaluated at $\widetilde B$.
Then
\begin{equation}\label{eq:coefficient-scaling}
 \dot G^{ij}(B)=b^{-(k+1)}\dot G^{ij}(\widetilde B),
 \qquad \mathcal C_G=b^{-k}\widetilde{\mathcal C_G}.
\end{equation}
With the derivatives of $G$ evaluated at $\widetilde B$,
\begin{align}
 (\partial_\tau-\dot G^{ij}\bar\nabla_i\bar\nabla_j)
       E_k(\widetilde B^{-1})
  =&(\dot G^{ij}\sigma_{ij}-k\overline F)E_k(\widetilde B^{-1}),
 \label{eq:normalized-speed-evolution}\\
 (\partial_\tau-\dot G^{ij}\bar\nabla_i\bar\nabla_j)
       E_1(\widetilde B)
 \notag=&-\frac2{E_1(\widetilde B)}\dot G^{ij}
       \bar\nabla_iE_1(\widetilde B)\bar\nabla_jE_1(\widetilde B)
       +(k-1)E_k(\widetilde B^{-1})
 \notag\\
 &+(\overline F-\dot G^{ij}\sigma_{ij})E_1(\widetilde B)
       -\frac1n\widetilde{\mathcal C_G},
 \label{eq:normalized-radius-evolution}
\end{align}
and
\begin{equation}\label{eq:normalized-product}
\begin{aligned}
 (\partial_\tau-\dot G^{ij}\bar\nabla_i\bar\nabla_j)
       \bigl(E_k(\widetilde B^{-1})E_1(\widetilde B)\bigr)=&-2\dot G^{ij}\bar\nabla_i\log E_1(\widetilde B)
        \bar\nabla_j\bigl(E_k(\widetilde B^{-1})E_1(\widetilde B)\bigr)\\
 &+(k-1)E_k(\widetilde B^{-1})^2
       -(k-1)\overline F E_k(\widetilde B^{-1})E_1(\widetilde B)
       \\&-\frac{E_k(\widetilde B^{-1})}{n}\widetilde{\mathcal C_G}.
\end{aligned}
\end{equation}
\end{lemma}

\begin{proof}
The relations \eqref{eq:scaling-identities} give
\eqref{eq:normalized-support}. Homogeneity gives
\begin{equation*}
 \dot G^{ij}(b\widetilde B)=b^{-(k+1)}\dot G^{ij}(\widetilde B),
 \qquad
 \ddot G^{ij,pq}(b\widetilde B)
 =b^{-(k+2)}\ddot G^{ij,pq}(\widetilde B),
\end{equation*}
and hence \eqref{eq:coefficient-scaling}.
Rescaling \eqref{eq:evolution-F} and
\eqref{eq:mean-radius-evolution} yields
\eqref{eq:normalized-speed-evolution} and
\eqref{eq:normalized-radius-evolution}.
The product rule proves \eqref{eq:normalized-product}.
\end{proof}

\begin{lemma}\label{lem:normalized-regularity}
Let $2\leq k\leq n-1$. Suppose that a solution of
\eqref{eq:normalized-support} satisfies
\begin{equation*}
 \|u(\cdot,\tau)\|_{C^0}\leq C,
 \qquad C^{-1}\sigma\leq\widetilde B(\cdot,\tau)\leq C\sigma
\end{equation*}
for $\tau_0\leq\tau<\tau_1$.
For all nonnegative integers $j,l$, there is a constant
$C_{j,l}=C_{j,l}(n,k,C,j,l)$ such that
\begin{equation}\label{eq:interior-regularity}
 \|\partial_\tau^j u(\cdot,\tau)\|_{C^l(\mathbb S^n)}
 \leq C_{j,l}\qquad(\tau_0+1\leq\tau<\tau_1).
\end{equation}
\end{lemma}

\begin{proof}
The radius bounds imply $C'^{-1}\leq\overline F\leq C'$.
On each interval $[a,a+1]$ contained in the given time interval,
set
\begin{equation*}
 \lambda(\tau)=\exp\left(-\int_a^\tau\overline F(\xi)\,d\xi\right),
 \qquad v=\lambda u,
 \qquad \frac{ds}{d\tau}=\lambda^{k+1}.
\end{equation*}
Then
\begin{equation*}
 \partial_sv=G(\bar\nabla^2v+v\sigma).
\end{equation*}
The spatial scaling and time change are uniformly bounded above
and below on this interval. Since $G$ is concave, the transformed
equation is uniformly parabolic and concave in the Hessian.
The parabolic Evans--Krylov estimate~\cite{Krylov1983},
followed by the parabolic Schauder estimates
\cite{Lieberman1996}*{Chapter~IV},
gives uniform bounds for all spatial derivatives away
from the initial time.  These transfer to $u$.
Differentiating \eqref{eq:normalized-support} and
\eqref{eq:average-normalized-speed} then gives the time derivative
bounds. Applying the estimates on preceding intervals of length
one proves \eqref{eq:interior-regularity}.
\end{proof}

\section{Sharp Codazzi inequalities for \texorpdfstring{$E_k$}{E k}}\label{sec:codazzi}

The Codazzi identity \eqref{eq:codazzi-radius} makes
$\bar\nabla B$ fully symmetric. We use this symmetry to control
both the derivatives of the mean principal radius and the gradient
of the speed by the Hessian of $G(B)=-E_k(B^{-1})$.
The use of fully symmetric third-order tensors to estimate gradient
terms in curvature flows is related to the pinching arguments in
\cite{Andrews2007}*{Sections~4--6} and
\cite{LuoWeiZhou}*{Lemma~2.1}.
The additional gradient term below provides a quantitative
lower bound for the rate of decrease of the monotone integral.

\begin{proposition}\label{prop:strict-codazzi}
Let $n\geq2$, $1\leq k<n$, and
$G(B)=-E_k(B^{-1})$ for $B>0$. Write $F=-G>0$.
For every fully symmetric three-tensor $T$,
\begin{equation}\label{eq:strengthened-codazzi}
 -\ddot G^{ij,pq}T_{mij}T_{mpq}
 \geq
 \frac{2}{nE_1(B)}\dot G^{ij}T_{ipp}T_{jqq}
 +\frac1F\sum_m\bigl(\dot G^{ij}T_{mij}\bigr)^2.
\end{equation}
Moreover, the weaker inequality
\begin{equation}\label{eq:strict-codazzi}
 -\ddot G^{ij,pq}T_{mij}T_{mpq}
 \geq
 \frac{2}{nE_1(B)}\dot G^{ij}T_{ipp}T_{jqq}
\end{equation}
has equality if and only if $T=0$.
On every compact subset $\mathcal K$ of the positive definite cone,
the difference between the two sides of \eqref{eq:strict-codazzi}
is bounded below by $c|T|^2$, where $c=c(n,k,\mathcal K)>0$.
The coefficient $2$ in \eqref{eq:strict-codazzi} is optimal. 
\end{proposition}

\begin{proof}
It suffices to work with $f(B)=\sigma_k(B^{-1})$, since
multiplication by a positive constant scales every term by the
same factor. Write $s=\operatorname{tr}B$ and $v_i=T_{ipp}$.
Repeated tensor indices are summed unless otherwise indicated.
Restricted sums are written explicitly.

\medskip
\noindent\emph{Step 1. The case $k=n-1$.}
Choose an orthonormal frame in which $B_{ij}=r_i\delta_{ij}$.
Then
\begin{equation*}
 f=\frac{s}{\det B},
 \qquad
 \log f=\log s-\log\det B,
 \qquad
 \dot f^{ii}=f\left(\frac1s-\frac1{r_i}\right)<0,
\end{equation*}
where there is no summation in the last identity.
Differentiating $\log f$ gives
\begin{equation}\label{eq:codazzi-basic-identity}
\begin{aligned}
 &\ddot f^{ij,pq}T_{mij}T_{mpq}
       +\frac2s\dot f^{ij}v_iv_j=
 \frac1f\sum_m\bigl(\dot f^{ij}T_{mij}\bigr)^2
 +f\left[
 \sum_{m,i,j}\frac{T_{mij}^2}{r_ir_j}
 -\sum_i\left(\frac2{sr_i}-\frac1{s^2}\right)v_i^2
 \right].
\end{aligned}
\end{equation}
The full symmetry of $T$ implies
\begin{equation*}
\begin{aligned}
 \sum_{m,i,j}\frac{T_{mij}^2}{r_ir_j}
 ={}&\sum_i\left[
 \frac{T_{iii}^2}{r_i^2}
 +\sum_{j\ne i}
 \left(\frac1{r_j^2}+\frac2{r_ir_j}\right)T_{ijj}^2
 \right]+2\sum_{i<j<l}
 \left(\frac1{r_ir_j}+\frac1{r_ir_l}+\frac1{r_jr_l}\right)
 T_{ijl}^2.
\end{aligned}
\end{equation*}
Fix $i$, with no summation over $i$ in the next three displays,
and set
\begin{equation*}
 \mu_i=\frac{\sum_{j\ne i}T_{ijj}}{s-r_i},
 \qquad v_i=T_{iii}+(s-r_i)\mu_i.
\end{equation*}
Since $r_j\leq s-r_i$ for $j\ne i$, weighted Cauchy's inequality
gives
\begin{equation*}
 \sum_{j\ne i}
 \left(\frac1{r_j^2}+\frac2{r_ir_j}\right)T_{ijj}^2
 \geq\frac{2s-r_i}{r_i(s-r_i)}
       \sum_{j\ne i}\frac{T_{ijj}^2}{r_j}
 \geq\frac{2s-r_i}{r_i}\mu_i^2.
\end{equation*}
Consequently,
\begin{equation}\label{eq:codazzi-row-square}
\begin{aligned}
 \frac{T_{iii}^2}{r_i^2}
 +\sum_{j\ne i}
 \left(\frac1{r_j^2}+\frac2{r_ir_j}\right)T_{ijj}^2\geq&\frac{T_{iii}^2}{r_i^2}
              +\frac{2s-r_i}{r_i}\mu_i^2\\
 =&
 \left(\frac2{sr_i}-\frac1{s^2}\right)v_i^2
 +\left(
 \frac{s-r_i}{sr_i}T_{iii}
 -\frac{2s-r_i}{s}\mu_i
 \right)^2.
\end{aligned}
\end{equation}
Summing over $i$ and using \eqref{eq:codazzi-basic-identity}
proves the strengthened inequality in this dimension.

If the weaker inequality is an equality, then
$\dot f^{pq}T_{ipq}=0$ for every $i$, all the squares in
\eqref{eq:codazzi-row-square} vanish, and equality holds in
weighted Cauchy's inequality.
Thus $T_{ijj}/r_j=\mu_i$ for $j\ne i$, and
\begin{equation*}
 \sum_j\frac{T_{ijj}}{r_j}-\frac{v_i}{s}=0,
 \qquad
 \frac{s-r_i}{sr_i}T_{iii}-\frac{2s-r_i}{s}\mu_i=0.
\end{equation*}
Subtracting the second expression from the first gives
$\mu_i=0$.
Hence $T_{iii}=T_{ijj}=0$ for all $i,j$.
The terms with three distinct indices also vanish, so $T=0$.

We will also use the following consequence of the same
logarithmic differentiation. For every symmetric two-tensor $V$,
\begin{equation}\label{eq:codazzi-log-convexity}
\begin{aligned}
 \ddot f^{ij,pq}V_{ij}V_{pq}
 -\frac1f\bigl(\dot f^{ij}V_{ij}\bigr)^2
 &=f\left[
 \sum_{i,j}\frac{V_{ij}^2}{r_ir_j}
 -\frac1{s^2}\left(\sum_iV_{ii}\right)^2
 \right]\\
 &\geq f\left(1-\frac{\sum_i r_i^2}{s^2}\right)
          \sum_{i,j}\frac{V_{ij}^2}{r_ir_j}.
\end{aligned}
\end{equation}
The last step is Cauchy's inequality, and the last expression
is strictly positive for $V\ne0$.

\medskip
\noindent\emph{Step 2. The general case.}
It remains to consider $k<n-1$.
Fix an orthonormal frame in which $B$ is diagonal.
For each $(k+1)$-element set
$U\subset\{1,\ldots,n\}$, let
$B_U=(B_{ij})_{i,j\in U}$ and set
\begin{equation*}
 s_U=\operatorname{tr}B_U,
 \qquad
 f_U=\sigma_k(B_U^{-1})=\frac{s_U}{\det B_U}.
\end{equation*}
For any positive definite matrix $B$, we have
\begin{equation*}
 \sigma_k(B^{-1})
 =\frac1{n-k}\sum_{|U|=k+1}
     \sigma_k\bigl((B^{-1})_U\bigr),
\end{equation*}
since each principal minor of order $k$ occurs
$n-k$ times in the sum.

At a diagonal matrix, $(B^{-1})_U=(B_U)^{-1}$,
and their first variations agree.
The derivatives of $f_U$ are taken with respect to
$B_U$ and extended by zero whenever an index lies
outside $U$.
In particular,
\begin{equation*}
 \dot f_U^{ii}
 =f_U\left(\frac1{s_U}-\frac1{r_i}\right)<0
 \qquad(i\in U),
\end{equation*}
with no summation over $i$.
The preceding identity therefore gives
\begin{equation}\label{eq:subset-counting}
 f=\frac1{n-k}\sum_U f_U,
 \qquad
 \dot f^{ij}=\frac1{n-k}\sum_U\dot f_U^{ij}.
\end{equation}

For the second variation, we must also account for
the indices outside $U$.
For any symmetric matrix $V$, we have
\begin{equation*}
 [(B+\varepsilon V)^{-1}]_{ij}
 =\frac{\delta_{ij}}{r_i}
 -\varepsilon\frac{V_{ij}}{r_ir_j}
 +\varepsilon^2\sum_l\frac{V_{il}V_{lj}}{r_ir_lr_j}
 +O(\varepsilon^3).
\end{equation*}
The corresponding expansion of
$(B_U+\varepsilon V_U)^{-1}$ contains only the terms
with $l\in U$.
Comparing the two expansions gives
\begin{equation}\label{eq:codazzi-subspace-hessian}
\begin{aligned}
 \ddot f^{ij,pq}V_{ij}V_{pq}
 =\frac1{n-k}\sum_U\left[
 \ddot f_U^{ij,pq}V_{ij}V_{pq}
 -2\sum_{i\in U}\sum_{j\notin U}
       \frac{\dot f_U^{ii}}{r_j}V_{ij}^2
 \right].
\end{aligned}
\end{equation}
All sums over $U$ are over $(k+1)$-element sets.

Apply this identity with $V_{ij}=T_{mij}$ and sum over $m$.
For $m\in U$, apply Step~1 to the fully symmetric tensor
restricted to $U$.
For $m\notin U$, apply \eqref{eq:codazzi-log-convexity}
to $(T_{mij})_{i,j\in U}$.
In the transverse sum, retain the terms with $m=j$ and use
$T_{jij}=T_{ijj}$.
It follows that
\begin{equation}\label{eq:codazzi-subspace-bound}
\begin{aligned}
 \ddot f^{ij,pq}T_{mij}T_{mpq}
 \geq\frac1{n-k}\sum_U\Bigg\{
 &-2\sum_{i\in U}\dot f_U^{ii}
 \left[
 \frac{\bigl(\sum_{j\in U}T_{ijj}\bigr)^2}{s_U}
 +\sum_{j\notin U}\frac{T_{ijj}^2}{r_j}
 \right]\\
 &+\frac1{f_U}\sum_m
       \bigl(\dot f_U^{ij}T_{mij}\bigr)^2
 \Bigg\}.
\end{aligned}
\end{equation}
For each $i\in U$, Cauchy's inequality gives
\begin{equation*}
 \frac{\bigl(\sum_{j\in U}T_{ijj}\bigr)^2}{s_U}
 +\sum_{j\notin U}\frac{T_{ijj}^2}{r_j}
 \geq
 \frac{\bigl(\sum_jT_{ijj}\bigr)^2}
      {s_U+\sum_{j\notin U}r_j}
 =\frac{v_i^2}{s}.
\end{equation*} 
Since $\dot f_U^{ii}<0$, \eqref{eq:subset-counting} gives
\begin{equation*}
 \ddot f^{ij,pq}T_{mij}T_{mpq}
 \geq-\frac2s\dot f^{ij}v_iv_j
 +\frac1{n-k}\sum_U\frac1{f_U}
       \sum_m\bigl(\dot f_U^{ij}T_{mij}\bigr)^2.
\end{equation*}
For each fixed $m$, another application of Cauchy's inequality
and \eqref{eq:subset-counting} yields
\begin{equation*}
 \frac1{n-k}\sum_U
 \frac{\bigl(\dot f_U^{ij}T_{mij}\bigr)^2}{f_U}
 \geq\frac1f\bigl(\dot f^{ij}T_{mij}\bigr)^2.
\end{equation*}
This proves the strengthened inequality for $f$.

For strictness of the weaker inequality, choose $T_{abc}\ne0$
and a set $U$ containing $b,c$.
If $a\in U$, Step~1 is strict for the restricted tensor.
If $a\notin U$, the term
$\ddot f_U^{ij,pq}T_{aij}T_{apq}$ is strictly positive by
\eqref{eq:codazzi-log-convexity} and can be retained in addition
to the trace terms when proving the weaker inequality.
Thus the weaker inequality is strict for $T\ne0$.
Continuity and compactness, after imposing $|T|=1$, give the
stated lower bound on compact subsets of the positive definite
cone.

Finally, take $B=\operatorname{diag}(L,1,\ldots,1)$ and
$T_{111}=1$, with all other components zero. Then
$s=L+n-1$, $v_1=1$, and $v_i=0$ for $i\ne1$. Since
\begin{equation*}
 f(B)=\binom{n-1}{k}+\frac1L\binom{n-1}{k-1},
\end{equation*}
direct differentiation gives
\begin{equation*}
 \frac{
 \ddot f^{ij,pq}T_{mij}T_{mpq}
 -\dfrac1f\sum_m\bigl(\dot f^{ij}T_{mij}\bigr)^2
 }{
 -s^{-1}\dot f^{ij}v_iv_j
 }
 =
 \frac{L+n-1}{L}
 \left(2-\frac{k}{(n-k)L+k}\right)
 \longrightarrow2
 \qquad(L\longrightarrow\infty).
\end{equation*}
Here we used $k<n$.
Thus the coefficient $2$ in
\eqref{eq:strengthened-codazzi} is optimal when the
coefficient of the last term is fixed at $1$.
For the weaker inequality, the corresponding quotient is
$2(L+n-1)/L\to2$, so the coefficient $2$ in
\eqref{eq:strict-codazzi} is also optimal.
\end{proof}

\begin{remark}
The equality assertion in Proposition~\ref{prop:strict-codazzi}
concerns \eqref{eq:strict-codazzi}, not the strengthened inequality
\eqref{eq:strengthened-codazzi}.
For example, when $n=2$, $k=1$, and $B=\operatorname{Id}$,
the symmetric tensor with
$T_{111}=3$ and $T_{122}=T_{212}=T_{221}=1$, all other components
being zero, attains equality in \eqref{eq:strengthened-codazzi}.
\end{remark}

\begin{remark}[Relation to inverse-concavity]
Set $\Phi(B)=\log E_k(B^{-1})$.
The inverse-concavity of $\log E_k$ in the sense of Andrews
is precisely the convexity of $\Phi$. See the proof of
\cite{Andrews2007}*{Theorem~2.1}.
In logarithmic form, \eqref{eq:strengthened-codazzi} reads
\begin{align*}
 &-\ddot G^{ij,pq}T_{mij}T_{mpq}
 -\frac1F\sum_m\bigl(\dot G^{ij}T_{mij}\bigr)^2\\
 =&
 F\ddot\Phi^{ij,pq}T_{mij}T_{mpq}
 \geq \frac{2}{nE_1(B)}
       \dot G^{ij}T_{ipp}T_{jqq}
 \geq0.
\end{align*}
Inverse-concavity gives the nonnegativity of the middle
expression by summing over $m$.
The full symmetry of $T$ yields the additional trace term,
so \eqref{eq:strengthened-codazzi} is a quantitative
Codazzi refinement of this traced inverse-concavity estimate.
\end{remark}
\begin{corollary}\label{cor:codazzi-remainder}
Let $B=\bar\nabla^2h+h\sigma>0$ and $1\leq k<n$.
The remainder in \eqref{eq:strict-remainder} satisfies
\begin{equation}\label{eq:remainder-speed-gradient}
 \mathcal C_G\geq
 \frac{|\bar\nabla E_k(B^{-1})|_\sigma^2}{E_k(B^{-1})}.
\end{equation}
Moreover,
\begin{equation}\label{eq:remainder-positivity}
 \mathcal C_G\geq0,
 \qquad
 \mathcal C_G=0\quad\Longleftrightarrow\quad\bar\nabla B=0
\end{equation}
at each point.
For every compact subset $\mathcal K$ of the positive definite cone,
there is $c=c(n,k,\mathcal K)>0$ such that
\begin{equation}\label{eq:remainder-coercivity}
 \mathcal C_G\geq c|\bar\nabla B|^2
 \qquad(B\in\mathcal K).
\end{equation}
The identity $\mathcal C_G\equiv0$ on $\mathbb S^n$ holds if and
only if the hypersurface is a round sphere.
The same assertions hold for $\widetilde B$ and
$\widetilde{\mathcal C_G}$, with $E_k(B^{-1})$ replaced by
$E_k(\widetilde B^{-1})$.
\end{corollary}

\begin{proof}
By the Codazzi identity, $T_{mij}=\bar\nabla_mB_{ij}$ is fully
symmetric, with
\begin{equation*}
 T_{mii}=n\bar\nabla_mE_1(B),
 \qquad
 \dot G^{ij}T_{mij}=-\bar\nabla_mE_k(B^{-1}).
\end{equation*}
The assertions follow from Proposition~\ref{prop:strict-codazzi}.
If $\mathcal C_G\equiv0$, then $\bar\nabla B=0$.
Every parallel symmetric two-tensor on the round sphere is a
constant multiple of $\sigma$, so $B=r\sigma$ for some $r>0$.
It follows that
\begin{equation*}
 \bar\nabla^2(h-r)+(h-r)\sigma=0,
\end{equation*}
and hence $h(z)=r+\langle c,z\rangle$ for some
$c\in\mathbb R^{n+1}$.
The converse follows directly from \eqref{eq:strict-remainder}.
The normalized assertions follow by applying the same argument
to $\widetilde B$.
\end{proof}

\section{A monotone entropy and geometric estimates}
\label{sec:integral}

We now assume $n\geq3$ and $2\leq k \leq n-1$ and construct a
scale-invariant entropy whose monotonicity follows from
Proposition~\ref{prop:strict-codazzi}.
Recall the exponents in \eqref{eq:exponents}:
\begin{equation*}
 p=\frac{n-1}{k-1},\qquad
 q=p-1=\frac{n-k}{k-1}.
\end{equation*}
For a smooth closed strictly convex hypersurface $M$, define
\begin{equation}\label{eq:monotone-integral}
\begin{aligned}
 \mathcal E(M)
 &=\frac1{|\mathbb S^n|}
     \int_M F^pE_1(B)^q\,d\mu
   =\fint_{\mathbb S^n}Z\,d\sigma,\\
 Z&=\frac{E_1(B)^qE_{n-k}(B)^p}{(\det B)^q}.
\end{aligned}
\end{equation}
Here $F=E_k(W)$.

In the following, we prove the monotonicity of $\mathcal E$ and derive an
integral estimate for the weighted Codazzi remainder.
The resulting uniform bound for $\mathcal E$ yields control of
the inradius and diameter, which will be used to obtain
estimates for the normalized flow
\eqref{eq:mean-width-normalization}.

\subsection{The entropy and its monotonicity}
The integral \eqref{eq:monotone-integral} has a projection
interpretation that explains its exponents. For background on integral-geometric formulas for convex bodies, see \cite{Schneider2014}*{Section~4.4}.
Set $j=n-k$, and equip the Grassmann manifold $\operatorname{Gr}(n+1,j+1)$ with its rotation-invariant
probability measure.
For $\Pi\in\operatorname{Gr}(n+1,j+1)$, let
$\Omega_\Pi=\Omega|\Pi$ be the orthogonal projection of $\Omega$,
with outward unit normal $\nu_\Pi$ and boundary area element
$d\mu_\Pi$.
For every continuous function $\Phi$ on $\mathbb S^n$,
\begin{equation}\label{eq:projection-measure}
 \fint_{\mathbb S^n}\Phi(z)E_j(B(z))\,d\sigma(z)
 =
 \frac1{|\mathbb S^j|}
 \int_{\operatorname{Gr}(n+1,j+1)}
 \int_{\partial\Omega_\Pi}
 \Phi(\nu_\Pi)\,d\mu_\Pi\,d\Pi.
\end{equation}

Indeed, the support function of $\Omega_\Pi$ is the restriction
of $h$ to $\mathbb S(\Pi)=\Pi\cap\mathbb S^n$.
At $z\in\mathbb S(\Pi)$, its radius tensor is the orthogonal
compression of $B(z)$ to $H=\Pi\cap z^\perp$, so its area density
is $\det(B(z)|_H)$.
The Cauchy--Binet formula and rotational invariance give
\begin{equation*}
 \int_{\operatorname{Gr}(z^\perp,j)}
 \det(B(z)|_H)\,dH=E_j(B(z)),
\end{equation*}
where $dH$ is the invariant probability measure.
Integrating over the pairs $(z,H)$ in either order proves
\eqref{eq:projection-measure}.

Since $F\det B=E_j(B)$ and $p=q+1$, we obtain
\begin{equation}\label{eq:projection-integral}
 \mathcal E(M)
 =
 \frac1{|\mathbb S^j|}
 \int_{\operatorname{Gr}(n+1,j+1)}
 \int_{\partial\Omega_\Pi}
 \bigl(FE_1(B)\bigr)^q(\nu_\Pi)\,d\mu_\Pi\,d\Pi.
\end{equation}
Here $F$ and $E_1(B)$ are evaluated on the original hypersurface
at the normal $\nu_\Pi$, rather than computed from
$\partial\Omega_\Pi$ alone. When $k=n-1$, the projections are planar convex bodies, and \eqref{eq:projection-integral} is the average of their weighted boundary lengths, divided by $2\pi$.

The evolution equation \eqref{eq:product-original}
suggests considering the product $FE_1(B)$.
Under dilation by $\lambda$, this product is multiplied by
$\lambda^{1-k}$, whereas $d\mu_\Pi$ is multiplied by $\lambda^j$.
Thus scale invariance of \eqref{eq:projection-integral} requires
$j-(k-1)q=0$, giving
$q=(n-k)/(k-1)$ and $p=q+1$.
This provides a geometric motivation for the integral. Its monotonicity uses the strict Codazzi inequality in Proposition \ref{prop:strict-codazzi}, as we
now prove.

\begin{proposition}\label{prop:monotonicity}
(i). The entropy $\mathcal E$ is invariant under translations and
dilations, and $\mathcal E\geq1$. Moreover, $\mathcal E(M)=1$ if and only if $M$ is a round sphere. 

(ii). Along a smooth closed strictly convex solution of \eqref{eq:flow},
write $w=FE_1(B)$. Then
\begin{equation}\label{eq:integral-derivative-original}
\begin{aligned}
 \frac d{dt}\mathcal E(M_t)
 ={}&-(n-k)\fint_{\mathbb S^n}F^2w^{q-1}
       \bigl(E_1(B)E_{n-k-1}(B)-E_{n-k}(B)\bigr)\,d\sigma\\
 &-pq\fint_{\mathbb S^n}w^{q-2}E_{n-k}(B)\dot G^{ij}
             \bar\nabla_iw\bar\nabla_jw\,d\sigma\\
 &-\frac qn\fint_{\mathbb S^n}Fw^{q-1}E_{n-k}(B)
             \mathcal C_G\,d\sigma\leq0.
\end{aligned}
\end{equation}
Equality holds at a time if and only if $M_t$ is a round sphere.
\end{proposition}

\begin{proof}
(i). Translation invariance is immediate, while dilation invariance follows
from $n-kp+q=0$. Since $(n-k)p=(n-1)q$, the
Newton--Maclaurin inequalities give
\begin{equation*}
Z=
\left(\frac{E_1(B)}{(\det B)^{1/n}}\right)^q
\left(\frac{E_{n-k}(B)}
{(\det B)^{(n-k)/n}}\right)^p
\geq1.
\end{equation*}
Hence $\mathcal E\geq1$. 
If $\mathcal E(M)=1$, then $Z=1$ everywhere.
Equality in the arithmetic-geometric mean inequality
gives $B=r\sigma$ for some positive function $r$.
The Codazzi equation implies that $r$ is constant,
so $M$ is a round sphere.
The converse is immediate.

(ii). Let $N^{ij}=\partial E_{n-k}(B)/\partial B_{ij}$.
By \eqref{eq:radius-newton-divergence} and
$\partial_tB=-\bar\nabla^2F-F\sigma$, we have
\begin{equation*}
\partial_tE_{n-k}(B)
=-\bar\nabla_i(N^{ij}\bar\nabla_jF)
-(n-k)FE_{n-k-1}(B).
\end{equation*}
On the other hand, differentiating
$E_{n-k}(B)=F\det B$ gives
\begin{equation*}
N^{ij}=(\det B)
\bigl(F(B^{-1})^{ij}-\dot G^{ij}\bigr).
\end{equation*}
Since both $N^{ij}$ and $(\det B)(B^{-1})^{ij}$
are divergence free, we obtain
\begin{equation}\label{eq:projection-divergence}
N^{ij}\bar\nabla_iF
=\bar\nabla_i\bigl(E_{n-k}(B)\dot G^{ij}\bigr)
-2(\det B)\dot G^{ij}\bar\nabla_iF.
\end{equation}

We now differentiate
$\mathcal E(M_t)=\fint_{\mathbb S^n}w^qE_{n-k}(B)\,d\sigma$.
Using \eqref{eq:product-original}, we obtain
\begin{align*}
\frac{d}{dt}\mathcal E(M_t)
=&q\fint_{\mathbb S^n}w^{q-1}E_{n-k}(B)
\bigl(
\dot G^{ij}\bar\nabla_i\bar\nabla_jw
-2\dot G^{ij}\bar\nabla_i\log E_1(B)\bar\nabla_jw
\bigr)\,d\sigma\\
&-\fint_{\mathbb S^n}w^q
\bar\nabla_i(N^{ij}\bar\nabla_jF)\,d\sigma\\
&+\fint_{\mathbb S^n} \left[q(k-1)
F^2w^{q-1}E_{n-k}(B)-(n-k)Fw^qE_{n-k-1}(B)\right]\,d\sigma\\
&-\frac qn\fint_{\mathbb S^n}
Fw^{q-1}E_{n-k}(B)\mathcal C_G\,d\sigma.
\end{align*}
Integrating by parts in the first two lines and using $q(k-1)=(n-k)$ yields
\begin{align*}
\frac{d}{dt}\mathcal E(M_t)
=&-q(q-1)\fint_{\mathbb S^n}
w^{q-2}E_{n-k}(B)\dot G^{ij}
\bar\nabla_iw\bar\nabla_jw\,d\sigma\\
&-q\fint_{\mathbb S^n}w^{q-1}
\Bigl[
\bar\nabla_i(E_{n-k}(B)\dot G^{ij})
-N^{ij}\bar\nabla_iF+2E_{n-k}(B)\dot G^{ij}
\bar\nabla_i\log E_1(B)
\Bigr]\bar\nabla_jw\,d\sigma\\
&-(n-k)\fint_{\mathbb S^n}
F^2w^{q-1}\left(E_1(B)E_{n-k-1}(B)-E_{n-k}(B)\right)\,d\sigma\\
&-\frac qn\fint_{\mathbb S^n}
Fw^{q-1}E_{n-k}(B)\mathcal C_G\,d\sigma.
\end{align*}
By \eqref{eq:projection-divergence} and
$w=FE_1(B)$, we have
\begin{align*}
&\bar\nabla_i(E_{n-k}(B)\dot G^{ij})
-N^{ij}\bar\nabla_iF
+2E_{n-k}(B)\dot G^{ij}\bar\nabla_i\log E_1(B)\\
=&2(\det B)\dot G^{ij}
\bigl(\bar\nabla_iF+F\bar\nabla_i\log E_1(B)\bigr)\\
=&2w^{-1}E_{n-k}(B)\dot G^{ij}\bar\nabla_iw.
\end{align*}
The two gradient contributions combine with coefficient
$-q(q-1)-2q=-q(q+1)=-pq$. Using $p=q+1$, we obtain
\eqref{eq:integral-derivative-original}.

Finally, the Newton--Maclaurin inequality
$E_1(B)E_{n-k-1}(B)\geq E_{n-k}(B)$, together with the positive definiteness of $(\dot G^{ij})$ and $\mathcal C_G\geq0$ from Corollary~\ref{cor:codazzi-remainder},
shows that $\mathcal E(M_t)$ is nonincreasing.
Equality holds only if $\mathcal C_G\equiv0$, which forces $M_t$ to be a round sphere by the same corollary. Conversely, equality holds for round spheres.
\end{proof}
For $2\leq k\leq n-1$, the equality case gives another proof of the
case $\alpha=1$ of the rigidity theorem of
Gao--Li--Ma~\cite{GaoLiMa2018}*{Theorem~1.4}.

\begin{corollary}\label{cor:shrinker-rigidity}
Let $n\geq3$ and $2\leq k\leq n-1$.
Every smooth closed strictly convex self-similar shrinking
solution of \eqref{eq:flow} is a family of round spheres.
\end{corollary}

\begin{proof}
Along a self-similar shrinking solution, the hypersurfaces
differ only by translations and dilations.
Thus $\mathcal E$ is constant by
Proposition~\ref{prop:monotonicity}.
The equality case of that proposition implies that
each hypersurface is a round sphere.
\end{proof}
For the normalization \eqref{eq:mean-width-normalization}, write
$\mathcal E(\tau)=\mathcal E(Y(\tau))$.
By scale invariance,
\begin{equation}\label{eq:integrand-scale-invariance}
 \mathcal E(\tau(t))=\mathcal E(M_t)
 =\fint_{\mathbb S^n}E_k(\widetilde B^{-1})^p
       E_1(\widetilde B)^q\det\widetilde B\,d\sigma.
\end{equation}
The density $Z$ is unchanged by rescaling.

\begin{corollary}\label{cor:normalized-integral}
Along the normalized flow \eqref{eq:mean-width-normalization}, let
\begin{equation}\label{eq:integrated-remainder}
\mathcal D(\tau)=\frac qn\fint_{\mathbb S^n}
E_k(\widetilde B^{-1})^pE_1(\widetilde B)^{q-1}
\widetilde{\mathcal C_G}\det\widetilde B\,d\sigma.
\end{equation}
Then
\begin{equation}\label{eq:normalized-derivative}
\frac{d}{d\tau}\mathcal E(\tau)\leq-\mathcal D(\tau)\leq0.
\end{equation}
In particular,
\begin{equation}\label{eq:integral-global-bound}
1\leq\mathcal E(\tau)\leq\mathcal E(M_0).
\end{equation}
\end{corollary}

\begin{proof}
Since $dt/d\tau=b^{k+1}$, the chain rule gives $\frac{d}{d\tau}\mathcal E(\tau)
=b^{k+1}\frac{d}{dt}\mathcal E(M_t)$.
Using \eqref{eq:scaling-identities},
\eqref{eq:coefficient-scaling}, and
$q(k-1)=n-k$, we obtain
\begin{equation*}
b^{k+1}Fw^{q-1}E_{n-k}(B)\mathcal C_G=E_k(\widetilde B^{-1})^p
E_1(\widetilde B)^{q-1}
\widetilde{\mathcal C_G}\det\widetilde B,
\end{equation*}
where $p=q+1$.
Thus \eqref{eq:integral-derivative-original} yields \eqref{eq:normalized-derivative}. The bounds in \eqref{eq:integral-global-bound}
follow from Proposition~\ref{prop:monotonicity}.
\end{proof}

When $k=n-1$, the first integral in
\eqref{eq:integral-derivative-original} vanishes identically.
The strict remainder in \eqref{eq:integrated-remainder} is
available throughout $2\leq k \leq n-1$.

\subsection{Inradius estimate}

We aim to control the inradius in terms of $\mathcal E$.

\begin{lemma}\label{lem:algebraic-coercivity}
For every $B>0$,
\begin{align}
 E_{n-k}(B)
 &\geq E_1(B)^{k/(n-1)}
       (\det B)^{(n-k-1)/(n-1)},
 \label{eq:symmetric-lower}\\
 Z&\geq
 \left(\frac{E_1(B)^n}{\det B}\right)^{1/(k-1)}.
 \label{eq:Z-determinant-control}
\end{align}
\end{lemma}

\begin{proof}
The Newton inequalities imply that $\log E_j(B)$ is
concave in $j$. Hence
\begin{equation*}
 \log E_{n-k}(B)
 \geq \frac{k}{n-1}\log E_1(B)
 +\frac{n-k-1}{n-1}\log\det B,
\end{equation*}
which proves \eqref{eq:symmetric-lower}. Therefore,
\begin{align*}
 Z &\geq E_1(B)^{q+kp/(n-1)}
       (\det B)^{-q+(n-k-1)p/(n-1)}\\
 &=\left(\frac{E_1(B)^n}{\det B}\right)^{1/(k-1)}.
\end{align*}
\end{proof}

The ratio $E_1(B)/(\det B)^{1/n}$ is the arithmetic mean
of the principal radii divided by their geometric mean.
Integrating \eqref{eq:Z-determinant-control} gives
\begin{equation*}
 \fint_{\mathbb S^n}
 \left[
 \left(\frac{E_1(B)}{(\det B)^{1/n}}\right)^{n/(k-1)}-1
 \right]d\sigma
 \leq\mathcal E(\partial\Omega)-1.
\end{equation*}
Thus the entropy deficit controls, in an integral sense,
the failure of the principal radii to be equal.

\begin{proposition}\label{prop:noncollapse}
Let $n\geq3$ and $2\leq k\leq n-1$. Let $\Omega\subset\mathbb R^{n+1}$ be a smooth
strictly convex body with support function $h$ and principal-radius
tensor $B=\bar\nabla^2h+h\sigma$.
Set
\begin{equation*}
 b=\fint_{\mathbb S^n}h\,d\sigma
  =\fint_{\mathbb S^n}E_1(B)\,d\sigma.
\end{equation*}
Then
\begin{equation}\label{eq:diam-mean-support}
 \operatorname{diam}\Omega\leq C(n)b.
\end{equation}
For each $e\in\mathbb S^n$, write
\begin{equation*}
 w_e=h(e)+h(-e),
 \qquad
 \psi(z)=\sqrt{1-\langle z,e\rangle^2}.
\end{equation*}
We have
\begin{equation}\label{eq:projection-width}
 \fint_{\mathbb S^n}\det B\,\psi\,d\sigma
 \leq C(n)w_e b^{n-1},
\end{equation}
and
\begin{equation}\label{eq:width-inradius}
 \min_{e\in\mathbb S^n}w_e
 \leq(n+2)r_-(\Omega).
\end{equation}
Moreover, there is a constant $c=c(n,k)>0$ such that
\begin{equation}\label{eq:noncollapse}
 r_-(\Omega)\geq cb\,\mathcal E(\partial\Omega)^{-(k-1)}.
\end{equation}
\end{proposition}

\begin{proof}
The equality defining $b$ follows by integrating
$\operatorname{tr}_{\sigma}B=\bar\Delta h+nh$.
Moreover, $b$ is invariant under translations.
After translating the midpoint of a diameter to the origin and
rotating, we may assume that $\Omega$ contains the segment joining
$-\frac12\operatorname{diam}\Omega\,e_1$ and
$\frac12\operatorname{diam}\Omega\,e_1$.
The definition of support function gives $h(z)\ge \frac{1}{2}\operatorname{diam}\Omega|z_1|$, hence we have
\begin{equation*}
 b\geq
 \frac12\operatorname{diam}\Omega
 \fint_{\mathbb S^n}|z_1|\,d\sigma,
\end{equation*}
which proves \eqref{eq:diam-mean-support}.

Fix $e\in\mathbb S^n$.
For a unit vector $v\in e^\perp$, the projection
$\Omega|v^\perp$ lies in a box with one side of length $w_e$
and the other $n-1$ sides of length at most
$\operatorname{diam}\Omega$.
Cauchy's projection formula
\cite{Schneider1998}*{Section~2.4, formula~(20)} yields
\begin{equation*}
 \int_{\partial\Omega}|\langle\nu,v\rangle|\,d\mu
 =2\operatorname{Vol}_n(\Omega|v^\perp)
 \leq2w_e(\operatorname{diam}\Omega)^{n-1}.
\end{equation*}
Averaging over $v\in\mathbb S^{n}\cap e^\perp$ and using rotational symmetry in $e^\perp$, we obtain
\begin{align*}
\fint_{\mathbb S^{n}\cap e^\perp}
 \int_{\partial\Omega}|\langle\nu,v\rangle|\,d\mu\,dv
=&
 \int_{\partial\Omega}\left(
 \fint_{\mathbb S^{n}\cap e^\perp}
 |\langle\nu,v\rangle|\,dv\right)d\mu \\
=&
 c_n\int_{\partial\Omega}
 \sqrt{1-\langle\nu,e\rangle^2}\,d\mu  \\
=&
 c_n\int_{\mathbb S^n}\det B\,\psi\,d\sigma,
\end{align*}
where $c_n>0$ depends only on $n$.
Consequently, \eqref{eq:diam-mean-support} gives
\begin{equation*}
 c_n\int_{\mathbb S^n}\det B\,\psi\,d\sigma
 \leq 2w_e(\operatorname{diam}\Omega)^{n-1}
 \leq C(n)w_eb^{n-1}.
\end{equation*}
This proves \eqref{eq:projection-width}.

The standard comparison between minimal width and inradius
\cite{AGN2016}*{(2.8)} implies \eqref{eq:width-inradius}.

Finally, by \eqref{eq:Z-determinant-control} and H\"older's
inequality with exponents $n/(k-1)$, $n$, and $n/(n-k)$,
\begin{equation*}
\begin{aligned}
 b^n
 &\leq C(n,k)\mathcal E^{k-1}
   \left(\fint_{\mathbb S^n}\det B\,\psi\,d\sigma\right)
   \left(\fint_{\mathbb S^n}
          \psi^{-1/(n-k)}\,d\sigma\right)^{n-k}\\
 &\leq C(n,k)\mathcal E^{k-1}w_e b^{n-1}.
\end{aligned}
\end{equation*}
Here we used \eqref{eq:projection-width}.
The integral involving $\psi$ is finite, since
$1/(n-k)\leq1<n$.
Thus
\begin{equation}\label{eq:width-lower}
 w_e\geq c(n,k)b\,\mathcal E^{-(k-1)}
 \qquad\text{for every }e\in\mathbb S^n.
\end{equation}
The inradius estimate now follows from
\eqref{eq:width-inradius}.
\end{proof}

\subsection{Quantitative stability of the entropy}
\label{subsec:monotone-stability}

The entropy also provides a quantitative measure of the
deviation of a convex body from a round ball.
The relevant geometric quantity is the mean area of its
two-dimensional projections.

For a smooth strictly convex body
$\Omega\subset\mathbb R^{n+1}$, let $h$ be its support
function and set
\begin{equation*}
 b=\fint_{\mathbb S^n}h\,d\sigma,
 \qquad
 c=(n+1)\fint_{\mathbb S^n}h(z)z\,d\sigma.
\end{equation*}
Thus $b$ is half the mean width and $c$ is the Steiner point.
Define the normalized mean area of the planar projections by
\begin{equation*}
 \mathcal A(\Omega)
 =\frac1\pi
 \int_{\operatorname{Gr}(n+1,2)}
       \operatorname{Area}(\Omega|\Pi)\,d\Pi,
\end{equation*}
where $d\Pi$ is the rotation-invariant probability measure.
For a ball of radius $b$, we have $\mathcal A(\Omega)=b^2$.

\begin{proposition}\label{prop:monotone-stability}
Let $n\geq3$ and $2\leq k\leq n-1$, and let
$\Omega\subset\mathbb R^{n+1}$ be a smooth strictly convex body.
Then
\begin{equation}\label{eq:stability-area}
 \mathcal E(\partial\Omega)^{-(k-1)}
 \leq
 \frac{|\partial\Omega|}{|\mathbb S^n|b^n}
 \leq1,
\end{equation}
and
\begin{equation}\label{eq:stability-projection-area}
 \mathcal E(\partial\Omega)^{-2(k-1)/n}
 \leq\frac{\mathcal A(\Omega)}{b^2}
 \leq1.
\end{equation}
Moreover,
\begin{equation}\label{eq:stability-hausdorff}
\begin{aligned}
 \frac{d_H(\Omega,\overline B_b(c))}{b}
 &\leq C(n,k)
 \bigl(\mathcal E(\partial\Omega)-1\bigr)^{2/(n+2)}.
\end{aligned}
\end{equation}
In particular, $\mathcal E(\partial\Omega)=1$ if and only if
$\Omega$ is a round ball.
\end{proposition}

\begin{proof}
Write $B=\bar\nabla^2h+h\sigma$.
The planar support-function area formula and rotational
invariance give
\begin{equation}\label{eq:mean-projection-identity}
\begin{aligned}
 \mathcal A(\Omega)
 &=\fint_{\mathbb S^n}
   \left(h^2-\frac1n|\bar\nabla h|^2\right)d\sigma\\
 &=\fint_{\mathbb S^n}hE_1(B)\,d\sigma
  =\fint_{\mathbb S^n}E_2(B)\,d\sigma.
\end{aligned}
\end{equation}
The last equality follows by integration by parts.
The Alexandrov--Fenchel inequalities
\cite{Schneider2014}*{Section~7.3} imply
\begin{equation*}
 \left(\fint_{\mathbb S^n}\det B\,d\sigma\right)^{1/n}
 \leq\mathcal A(\Omega)^{1/2}\leq b.
\end{equation*}
On the other hand, \eqref{eq:Z-determinant-control}
and H\"older's inequality yield
\begin{equation*}
 b\leq\mathcal E(\partial\Omega)^{(k-1)/n}
 \left(\fint_{\mathbb S^n}\det B\,d\sigma\right)^{1/n}.
\end{equation*}
Combining these inequalities proves
\eqref{eq:stability-area} and
\eqref{eq:stability-projection-area}.

By translation and dilation, we may assume $c=0$ and $b=1$.
Write $h=1+v$. Then $v$ is orthogonal to the constant and
first spherical harmonics.
By \eqref{eq:mean-projection-identity} and the spectral gap
of $-\bar\Delta$ on their orthogonal complement,
\begin{equation}\label{eq:stability-H1}
 \|v\|_{H^1(\mathbb S^n)}^2
 \leq C(n)\bigl(1-\mathcal A(\Omega)\bigr)
 \leq C(n)
 \left[1-\mathcal E(\partial\Omega)^{-2(k-1)/n}\right].
\end{equation}
All norms use the normalized spherical measure.

The inverse Gauss map and integration by parts give
\begin{equation*}
 \fint_{\mathbb S^n}X\,d\sigma
 =\fint_{\mathbb S^n}(hz+\bar\nabla h)\,d\sigma
 =c=0,
\end{equation*}
so $0\in\Omega$.
Set $R=\max_{\mathbb S^n}h$ and
$a=\|h-1\|_{C^0(\mathbb S^n)}$.
By Proposition~\ref{prop:noncollapse},
\begin{equation*}
 0\leq h\leq R\leq C(n),
 \qquad \Omega\subset\overline B_R(0).
\end{equation*}
Since $X=hz+\bar\nabla h$ and $R-h\leq2a$,
\begin{equation*}
\begin{aligned}
 |\bar\nabla h|^2
 &=|X|^2-h^2\\
 &\leq(R-h)(R+h)
 \leq4Ra
 \leq C(n)a.
\end{aligned}
\end{equation*}
If $a>0$, integration along geodesics from a maximum point
of $|v|$ shows that $|v|\geq a/2$ on a spherical cap of
radius $c(n)\sqrt a$.
The measure of this cap is bounded below by $c(n)a^{n/2}$.
Since $n\geq3$, Sobolev's inequality gives
\begin{equation*}
 c(n)a^{(n+2)/4}
 \leq\|v\|_{L^{2n/(n-2)}(\mathbb S^n)}
 \leq C(n)\|v\|_{H^1(\mathbb S^n)}.
\end{equation*}
The case $a=0$ is immediate. 
Combining this with \eqref{eq:stability-H1} and
\eqref{eq:stability-projection-area}, we obtain
\begin{align*}
 \|h-1\|_{C^0(\mathbb S^n)}
 &\leq C(n)\bigl(1-\mathcal A(\Omega)\bigr)^{2/(n+2)}\\
 &\leq C(n)
 \left[1-\mathcal E(\partial\Omega)^{-2(k-1)/n}\right]^{2/(n+2)}\\
 &\leq C(n,k)
 \bigl(\mathcal E(\partial\Omega)-1\bigr)^{2/(n+2)}.
\end{align*}
Here we used $1-s^{-\beta}\leq\beta(s-1)$ for
$s\geq1$ and $\beta>0$.
Restoring the center and scale gives
\eqref{eq:stability-hausdorff}, since the Hausdorff
distance between convex bodies equals the uniform distance
between their support functions.
The equality statement follows as well.
\end{proof}

\begin{remark}
This proposition gives a static geometric interpretation
of the entropy: its deficit controls the Hausdorff distance
to the ball with the same mean width.
We record this quantitative stability as an independent
geometric consequence of the entropy.
It will not be used in the subsequent proof of convergence
of the normalized flow.
\end{remark}

\section{Smooth convergence to a round point}\label{sec:convergence}

We now work with the normalized flow
\eqref{eq:mean-width-normalization}, for which the mean of the
support function is one. Unless otherwise stated, all constants depend only on $M_0$, $n$, and $k$.

\begin{lemma}\label{lem:normalized-speed}
Let $2\leq k\leq n-1$. The normalized flow \eqref{eq:mean-width-normalization} is defined for all $\tau\geq0$ and satisfies
\begin{equation}\label{eq:normalized-geometry}
 r_-(\widetilde\Omega_\tau)\geq c,
 \qquad \operatorname{diam}Y(\tau)\leq C,\qquad 0<u\leq C,
\end{equation}
and
\begin{equation}\label{eq:speed-product-bounds}
 E_k(\widetilde B^{-1})\leq C,
 \qquad 1\leq\overline F\leq C,
 \qquad E_k(\widetilde B^{-1})E_1(\widetilde B)\leq C.
\end{equation}
Moreover,
\begin{equation}\label{eq:mean-width-extinction-comparison}
 c\ell(t)\leq b(t)\leq C\ell(t)\qquad(0\leq t<T).
\end{equation}
\end{lemma}

\begin{proof}
Since $\fint_{\mathbb S^n}u\,d\sigma=1$,
Corollary~\ref{cor:normalized-integral} and
Proposition~\ref{prop:noncollapse} imply
\eqref{eq:normalized-geometry}.
By strict nesting and extinction at $x_*$, we have
$x_*\in\operatorname{int}\Omega_t$ for $t<T$.
Consequently, $0<u\leq\operatorname{diam}Y\leq C$. Notice that the inradius bound does not yield a uniform
positive lower bound for $u$, since $x_*$ need not be an
incenter.

Scaling back and using
\eqref{eq:extinction-radius-comparison}, we obtain
\begin{equation*}
cb(t)\leq r_-(\Omega_t)\leq\ell(t)
\leq\operatorname{diam}\Omega_t\leq Cb(t),
\end{equation*}
which proves \eqref{eq:mean-width-extinction-comparison}.
Since $\ell(t)^{k+1}$ is proportional to $T-t$, it follows that
\begin{equation*}
\frac{d\tau}{dt}=b(t)^{-(k+1)}
\geq\frac{c}{T-t}.
\end{equation*}
Hence $\tau(t)\to\infty$ as $t\to T$.

By Lemma~\ref{lem:harnack}, integrating the Harnack inequality
up to the extinction time as in
\cite{HamiltonEntropy1994}*{Section~2}, we obtain
\begin{align*}
h_{x_*}(z,t)
&=\int_t^T F(z,s)\,ds\geq F(z,t)\int_t^T
\left(\frac ts\right)^{k/(k+1)}ds\\
&\geq \left(\frac tT\right)^{k/(k+1)}
(T-t)F(z,t),
\end{align*}
where $h_{x_*}=h-\langle x_*,z\rangle$.
Consequently,
\begin{equation}\label{eq:normalized-speed-upper}
E_k(\widetilde B^{-1})
\leq
\left(\frac Tt\right)^{k/(k+1)}
\frac{b(t)^{k+1}}{T-t}\,u
\leq C
\end{equation}
for $T/2\leq t<T$, where we used
\eqref{eq:mean-width-extinction-comparison} and the
uniform upper bound for $u$.
By smoothness, the estimate also holds on $[0,T/2]$
after increasing $C$.
Together with \eqref{eq:average-speed-lower}, this yields $1\leq\overline F\leq C$.

Recall that
$w=E_k(\widetilde B^{-1})E_1(\widetilde B)$.
At a spatial maximum of $w$,
\eqref{eq:normalized-product} gives
\begin{equation*}
\partial_\tau w
\leq (k-1)E_k(\widetilde B^{-1})^2
-(k-1)\overline Fw  \leq (k-1)(C-w).
\end{equation*}
The maximum principle therefore implies $w\leq C$.
\end{proof}
\begin{lemma}\label{lem:support-from-dissipation}
Let $2\leq k \leq n-1$, and let $u$ be the normalized support function
in \eqref{eq:normalized-support-definition}. There is a constant
$C=C(M_0,n,k)>0$ such that, for every $\tau\geq0$,
\begin{equation}\label{eq:support-dissipation-tail}
 \fint_{\mathbb S^n}|u(z,\tau)-1|^2\,d\sigma
 \leq C\int_\tau^\infty e^{-(s-\tau)}\mathcal D(s)\,ds,
\end{equation}
where $\mathcal D$ is defined in \eqref{eq:integrated-remainder}.
In particular,
\begin{equation*}
 u(\cdot,\tau)\longrightarrow1
 \quad\text{uniformly on }\mathbb S^n
 \quad\text{as }\tau\to\infty.
\end{equation*}
\end{lemma}

\begin{proof}
In this proof, write $F=E_k(\widetilde{B}^{-1})$. Thus
$\overline F=\fint_{\mathbb S^n}F\,d\sigma$.
Lemma~\ref{lem:normalized-speed} gives
\begin{equation*}
 0<u\leq C,\qquad 0<F\leq C,\qquad
 1\leq\overline F\leq C,\qquad FE_1(\widetilde{B})\leq C.
\end{equation*}
Since $Z=F^pE_1(\widetilde{B})^q\det \widetilde{B}\geq1$, the strengthened inequality
\eqref{eq:remainder-speed-gradient} and the definition of
$\mathcal D$ imply
\begin{equation*}
 \mathcal D(\tau)
 \geq\frac qn\fint_{\mathbb S^n}
       \frac{Z}{FE_1(\widetilde{B})}|\bar\nabla F|^2\,d\sigma
 \geq c\fint_{\mathbb S^n}|\bar\nabla F|^2\,d\sigma.
\end{equation*}
The scalar Poincar\'e inequality therefore yields
\begin{equation*}
 \fint_{\mathbb S^n}|F-\overline F|^2\,d\sigma
 \leq C\mathcal D(\tau).
\end{equation*}

By \eqref{eq:normalized-support},
\begin{equation*}
 \partial_\tau(u-1)=\overline F(u-1)-(F-\overline F).
\end{equation*}
For $S>\tau$, integration gives
\begin{align*}
 u(z,\tau)-1
 =&e^{-\int_\tau^S\overline F(\xi)\,d\xi}
          \bigl(u(z,S)-1\bigr)\\
 &+\int_\tau^S e^{-\int_\tau^s\overline F(\xi)\,d\xi}
          \bigl(F(z,s)-\overline F(s)\bigr)\,ds.
\end{align*}
Since $u$ is uniformly bounded and $\overline F\geq1$,
the first term tends uniformly to zero as $S\to\infty$.
The integral converges absolutely because $F$ and $\overline F$
are bounded. Hence
\begin{equation*}
 u(z,\tau)-1
 =\int_\tau^\infty
       e^{-\int_\tau^s\overline F(\xi)\,d\xi}
       \bigl(F(z,s)-\overline F(s)\bigr)\,ds.
\end{equation*}
Using $\overline F\geq1$, we have
\begin{equation*}
 \int_\tau^\infty
       e^{-\int_\tau^s\overline F(\xi)\,d\xi}\,ds
 \leq\int_\tau^\infty e^{-(s-\tau)}\,ds=1.
\end{equation*}
The weighted Cauchy inequality and the speed estimate above give
\begin{align*}
 \fint_{\mathbb S^n}|u(z,\tau)-1|^2\,d\sigma
 &\leq\int_\tau^\infty
       e^{-\int_\tau^s\overline F(\xi)\,d\xi}
       \fint_{\mathbb S^n}
       |F(z,s)-\overline F(s)|^2\,d\sigma\,ds\\
 &\leq C\int_\tau^\infty e^{-(s-\tau)}\mathcal D(s)\,ds.
\end{align*}
This proves \eqref{eq:support-dissipation-tail}.

By \eqref{eq:normalized-derivative} and $\mathcal E\geq1$,
\begin{equation*}
 \int_0^\infty\mathcal D(s)\,ds
 \leq\mathcal E(0)-1<\infty.
\end{equation*}
Thus the right-hand side of \eqref{eq:support-dissipation-tail}
tends to zero, which proves the $L^2$ convergence.

Finally, since $0\in\widetilde\Omega_\tau$ and
$\operatorname{diam}\widetilde\Omega_\tau\leq C$,
\begin{equation*}
 |\bar\nabla u|^2+u^2
 =|\bar\nabla u+uz|^2\leq C.
\end{equation*}
Integrating over a spherical cap about a maximum point
of $|u-1|$ gives
\begin{align*}
 \|u(\cdot,\tau)-1\|_{C^0(\mathbb S^n)}^{n+2}
 &\leq C\fint_{\mathbb S^n}|u(z,\tau)-1|^2\,d\sigma\\
 &\leq C\int_\tau^\infty
             e^{-(s-\tau)}\mathcal D(s)\,ds\\
 &\leq C\int_\tau^\infty\mathcal D(s)\,ds
 \longrightarrow0.
\end{align*}
This proves uniform convergence.
\end{proof}

We now prove Theorem~\ref{thm:main}.

\begin{proof}[Proof of Theorem~\ref{thm:main}]
Point extinction follows from Proposition~\ref{prop:point-extinction}.
We use \eqref{eq:mean-width-normalization}. Lemma~\ref{lem:normalized-speed} ensures that $\tau\in[0,\infty)$. By Lemma~\ref{lem:support-from-dissipation},
there is $\tau_1$ such that $u\geq1/2$
for $\tau\geq\tau_1$. 
We use the test function from the alternative proof of
\cite{AGN2016}*{Theorem~5.2}.
For this step, write $F=E_k(\widetilde B^{-1})$, with derivatives
of $G$ evaluated at $\widetilde B$.
Equations \eqref{eq:normalized-support} and
\eqref{eq:normalized-speed-evolution} give
\begin{align*}
 (\partial_\tau-\dot G^{ij}\bar\nabla_i\bar\nabla_j)u
 &=(\overline F+\dot G^{ij}\sigma_{ij})u-(k+1)F,\\
 (\partial_\tau-\dot G^{ij}\bar\nabla_i\bar\nabla_j)F
 &=(\dot G^{ij}\sigma_{ij}-k\overline F)F.
\end{align*}
The gradient terms in the logarithmic equations are nonnegative.
Using $\overline F\geq1$, we obtain
\begin{equation}\label{eq:logarithmic-speed-test}
\begin{aligned}
 &(\partial_\tau-\dot G^{ij}\bar\nabla_i\bar\nabla_j)
                  \log(Fu^{k+1})
 \geq\overline F+(k+2)\dot G^{ij}\sigma_{ij}
                   -(k+1)^2\frac Fu
 \geq1-(k+1)^2\frac Fu.
\end{aligned}
\end{equation}
Let $m(\tau)=\min_{\mathbb S^n}Fu^{k+1}$.
At a spatial minimum,
$F/u=m(\tau)/u^{k+2}\leq2^{k+2}m(\tau)$.
The maximum principle gives
\begin{equation*}
 m(\tau)\geq\min\left\{
 m(\tau_1),\frac1{2^{k+3}(k+1)^2}\right\}>0
 \qquad(\tau\geq\tau_1).
\end{equation*}
Since $u\leq C$,
\begin{equation}\label{eq:normalized-speed-lower}
 E_k(\widetilde B^{-1})\geq c>0\qquad(\tau\geq\tau_1).
\end{equation}
The product bound in \eqref{eq:speed-product-bounds} now gives
$E_1(\widetilde B)\leq C$, and hence $\kappa_i(Y)\geq c$.
Keeping the terms containing the largest principal curvature,
\begin{equation*}
 C\geq E_k(\widetilde B^{-1})
 \geq\frac kn c^{k-1}\kappa_{\max}(Y).
\end{equation*}
Thus
\begin{equation}\label{eq:uniform-curvature}
 c\leq\kappa_i(Y)\leq C
\end{equation}
for all sufficiently large normalized times.
Lemma~\ref{lem:normalized-regularity} gives
\begin{equation}\label{eq:global-regularity}
 \|u(\cdot,\tau)\|_{C^l(\mathbb S^n)}\leq C_l
 \qquad(l\geq0).
\end{equation}
Interpolation with the $L^2$ convergence in
Lemma~\ref{lem:support-from-dissipation} proves smooth
convergence to one.

Finally, smooth convergence gives $\overline F(\tau)\to1$.
Since
\begin{equation*}
 \frac d{dt}b(t)^{k+1}=-(k+1)\overline F(\tau(t)),
 \qquad b(t)\longrightarrow0\quad(t\to T),
\end{equation*}
integration yields
\begin{equation}\label{eq:normalization-asymptotic}
 \frac{b(t)^{k+1}}{(k+1)(T-t)}\longrightarrow1.
\end{equation}
Using $\kappa_i(M_t)=b(t)^{-1}\kappa_i(Y(\tau(t)))$ and
$E_k(W)=b(t)^{-k}E_k(\widetilde W)$, we obtain
\eqref{eq:main-curvature} from
\eqref{eq:uniform-curvature}.
Increasing the upper constant covers the remaining compact
interval of original time.
\end{proof}

We finally prove that the convergence in \eqref{eq:main-convergence} is exponential in the normalized time.

\begin{corollary}\label{cor:exponential-convergence}
Let $2\leq k\leq n-1$, and let $u$ be the support function
in the mean-width normalization
\eqref{eq:mean-width-normalization}.
Set
\begin{equation*}
 \mu_*=k-1+\frac{2k}{n}.
\end{equation*}
For every $0<\lambda<\mu_*$ and every integer $j\geq0$, there is
$C=C(M_0,n,k,j,\lambda)>0$ such that
\begin{equation}\label{eq:exponential-convergence}
 \|u(\cdot,\tau)-1\|_{C^j(\mathbb S^n)}
 \leq Ce^{-\lambda\tau}\qquad(\tau\geq0).
\end{equation}
\end{corollary}

\begin{proof}
Write $F=E_k(\widetilde B^{-1})$ and
$\overline F=\fint_{\mathbb S^n}F\,d\sigma$.
By \eqref{eq:normalized-support} and
\eqref{eq:normalized-mean-support},
\begin{equation*}
 u_\tau=\overline F u-F,\qquad
 \fint_{\mathbb S^n}u\,d\sigma=1.
\end{equation*}
The smooth convergence proved above gives $u\to1$ in $C^\infty(\mathbb{S}^n)$.
In particular, all higher-order norms of $u$ are uniformly bounded.

Let $\Pi_m$ be the $L^2$ projection onto spherical harmonics of
degree $m$, and put
\begin{equation*}
 \eta=\Pi_1(u-1),\qquad
 w=(I-\Pi_0-\Pi_1)(u-1).
\end{equation*}
The mean-width constraint gives $u=1+\eta+w$. Since
$\bar\nabla^2\eta+\eta\sigma=0$, we have
\begin{equation*}
 \widetilde B=\sigma+Q,\qquad
 Q=\bar\nabla^2w+w\sigma.
\end{equation*}
Thus the curvature speed depends only on $w$.

Recall that $G(B)=-E_k(B^{-1})$ and
$\dot G^{ij}(\sigma)=\frac{k}{n}\sigma^{ij}$.
Taylor's formula gives
\begin{equation}\label{eq:exp-speed-expansion}
 F=1-\frac{k}{n}(\bar\Delta+n)w+R,
 \qquad \overline F=1+\overline R,
\end{equation}
where $\overline R=\fint_{\mathbb S^n}R\,d\sigma$ and
\begin{equation*}
 R=a^{ij}Q_{ij},\qquad
 a^{ij}=-\int_0^1\left(\dot G^{ij}(\sigma+sQ)
                  -\frac{k}{n}\sigma^{ij}\right)\,ds.
\end{equation*}
For all sufficiently large $\tau$, smoothness of $G$ and
$w\to0$ in $C^\infty$ imply
\begin{equation}\label{eq:exp-remainder-bounds}
 \|a\|_{C^1(\mathbb{S}^n)}\leq C\|w\|_{C^3(\mathbb{S}^n)}=o(1),
 \qquad
 \|R\|_{C^0(\mathbb{S}^n)}\leq C\|w\|_{C^2(\mathbb{S}^n)}^2.
\end{equation}

Projecting the support equation onto harmonics of degree at least
two, we obtain
\begin{equation*}
 w_\tau=\frac{k}{n}\bar\Delta w+(k+1)w
            +\overline R\,w-\Pi_{\geq2}R,
 \qquad \Pi_{\geq2}=I-\Pi_0-\Pi_1.
\end{equation*}
Taking the $L^2$ inner product with $w$ gives
\begin{equation*}
 \frac12\frac{d}{d\tau}\|w\|_{L^2(\mathbb{S}^n)}^2
 =-\frac{k}{n}\|\bar\nabla w\|_{L^2(\mathbb{S}^n)}^2
  +(k+1+\overline R)\|w\|_{L^2(\mathbb{S}^n)}^2
  -\int_{\mathbb S^n}wR\,d\sigma.
\end{equation*}
Integration by parts and \eqref{eq:exp-remainder-bounds} yield
\begin{equation*}
 \left|\int_{\mathbb S^n}wR\,d\sigma\right|
 \leq C\|a\|_{C^1(\mathbb{S}^n)}
       \bigl(\|\bar\nabla w\|_{L^2(\mathbb{S}^n)}^2+\|w\|_{L^2(\mathbb{S}^n)}^2\bigr).
\end{equation*}
Since $\overline R=o(1)$, for every $\varepsilon>0$ and all
sufficiently large $\tau$,
\begin{equation*}
 \frac12\frac{d}{d\tau}\|w\|_{L^2(\mathbb{S}^n)}^2
 \leq-\left(\frac{k}{n}-\varepsilon\right)
              \|\bar\nabla w\|_{L^2(\mathbb{S}^n)}^2
       +(k+1+\varepsilon)\|w\|_{L^2(\mathbb{S}^n)}^2.
\end{equation*}
The first eigenvalue of $-\bar\Delta$ on the space of spherical
harmonics of degree at least two is $2(n+1)$. Hence
\begin{equation*}
 \|\bar\nabla w\|_{L^2(\mathbb{S}^n)}^2
 \geq2(n+1)\|w\|_{L^2(\mathbb{S}^n)}^2,
\end{equation*}
and consequently
\begin{equation*}
 \frac12\frac{d}{d\tau}\|w\|_{L^2(\mathbb{S}^n)}^2
 \leq-\bigl(\mu_*-(2n+3)\varepsilon\bigr)
                  \|w\|_{L^2(\mathbb{S}^n)}^2.
\end{equation*}
It follows that, for every $0<\nu<\mu_*$,
\begin{equation}\label{eq:exp-high-modes}
 \|w(\cdot,\tau)\|_{L^2(\mathbb{S}^n)}\leq C_\nu e^{-\nu\tau}.
\end{equation}

Fix $0<\lambda<\mu_*$ and choose $\lambda<\nu<\mu_*$.
For any integer $j\geq0$, choose integers $N>s>j+n/2$ with
$\nu(1-s/N)>\lambda$. Sobolev embedding, interpolation, and the
uniform $H^N$ bound give
\begin{equation}\label{eq:exp-high-modes-Cj}
 \|w\|_{C^j(\mathbb{S}^n)}
 \leq C\|w\|_{H^s(\mathbb{S}^n)}
 \leq C\|w\|_{L^2(\mathbb{S}^n)}^{1-s/N}\|w\|_{H^N(\mathbb{S}^n)}^{s/N}
 \leq C_{j,\lambda}e^{-\lambda\tau}.
\end{equation}

It remains to control the first harmonics. Projecting
\eqref{eq:normalized-support} and using
\eqref{eq:exp-speed-expansion}, we obtain
\begin{equation*}
 \eta_\tau=\overline F\eta-\Pi_1R.
\end{equation*}
For $S>\tau$, variation of constants gives
\begin{equation*}
 \eta(\tau)
 =e^{-\int_\tau^S\overline F(\xi)\,d\xi}\eta(S)
  +\int_\tau^S e^{-\int_\tau^s\overline F(\xi)\,d\xi}
               \Pi_1R(s)\,ds.
\end{equation*}
Since $\overline F\geq1$ by \eqref{eq:average-speed-lower}
and $\eta$ is bounded, the first term tends to zero as
$S\to\infty$. Thus
\begin{equation*}
 \eta(\tau)=\int_\tau^\infty
       e^{-\int_\tau^s\overline F(\xi)\,d\xi}
       \Pi_1R(s)\,ds.
\end{equation*}
The first-harmonic space is finite dimensional. By
\eqref{eq:exp-remainder-bounds} and
\eqref{eq:exp-high-modes-Cj},
\begin{equation*}
 \|\Pi_1R(s)\|_{C^j(\mathbb{S}^n)}
 \leq C_j\|R(s)\|_{C^0(\mathbb{S}^n)}
 \leq C_j\|w(s)\|_{C^2(\mathbb{S}^n)}^2
 \leq C_{j,\lambda}e^{-2\lambda s}.
\end{equation*}
It follows that
\begin{equation*}
 \|\eta(\tau)\|_{C^j(\mathbb{S}^n)}
 \leq C_{j,\lambda}\int_\tau^\infty
                e^{-(s-\tau)}e^{-2\lambda s}\,ds
 \leq C_{j,\lambda}e^{-2\lambda\tau}.
\end{equation*}
Combining this with \eqref{eq:exp-high-modes-Cj} proves
\eqref{eq:exponential-convergence} for all sufficiently large
$\tau$. Enlarging the constant covers the remaining compact
time interval.
\end{proof}

\section*{Acknowledgements}

This work was supported by the National Key Research and Development Program of China (2021YFA1001800), the National Natural Science Foundation of China (Grant No.~12531002), and the Fundamental Research Funds for the Central Universities. The third author was also supported by the National Natural Science Foundation of China (Grant No.~12601099). The last author was also supported by the China Postdoctoral Science Foundation (Grant No.~2025M783146).

\end{document}